\documentclass[a4paper,12pt]{amsart}
\usepackage{amsfonts,graphicx}
\usepackage{subcaption}
\usepackage{amsmath, amssymb, amsthm,bm}
\usepackage{mathrsfs,enumitem} 
\usepackage{graphicx}
\usepackage{color}
\usepackage[hyperindex,breaklinks]{hyperref}
\usepackage{xspace}
\usepackage{framed}
\usepackage[round]{natbib}
\usepackage[noend]{algpseudocode}
\usepackage[makeroom]{cancel}
\usepackage[a4paper]{geometry}
\theoremstyle{definition}
\newtheorem{theorem}{Theorem}

\newtheorem{assumption}{Assumption}

\newtheorem{definition}{Definition}

\newtheorem{lemma}[theorem]{Lemma}

\newtheorem{proposition}[theorem]{Proposition}
\newtheorem{remark}{Remark}

\graphicspath{{.}{./figures/}}

\newenvironment{revisionenv}
{}
{}

\newenvironment{proofenv}[1]
{\begin{proof}[#1]}
{\end{proof}
}

\usepackage{caption}
\usepackage[font=small,position=top]{subcaption}
\usepackage{xcolor}
\usepackage{calrsfs}
\usepackage{hyperref}
\usepackage{enumitem}
 \usepackage[ruled,linesnumbered,vlined]{algorithm2e}
 \newtheorem{property}{Property}
 \DeclareMathOperator{\Val}{Val}
 \newcommand{\dconverge}{\overset{d}{\rightarrow}}
 \newcommand{\vconverge}{\overset{v}{\rightarrow}}

\usepackage{natbib}
 \bibpunct[, ]{(}{)}{,}{a}{}{,}%
\begin{document}

\title[Efficient Scenario Generation for Chance Constrained Optimization]{Efficient Scenario Generation for Heavy-tailed Chance Constrained Optimization}
	\author{Jose Blanchet}
	\address{Department of Management Science and Engineering, Stanford University, Stanford, CA 94305}
	\email{jose.blanchet@stanford.edu}
	\author{Fan Zhang}
	\address{Department of Management Science and Engineering, Stanford University, Stanford, CA 94305}
	\email{fzh@stanford.edu}
	\author{Bert Zwart}
	\address{Centrum Wiskunde \& Informatica, Amsterdam 1098 XG, Netherlands, and Eindhoven University of Technology.}
	\email{bert.zwart@cwi.nl}
	\maketitle
	\date{\today }
	
	\begin{abstract}
		We consider a generic class of chance-constrained optimization problems with heavy-tailed (i.e., power-law type) risk factors. In this setting, we use the scenario approach to obtain a constant approximation to the optimal solution with a computational complexity that is uniform in the risk tolerance parameter. We additionally illustrate the
    efficiency of our algorithm in the context of solvency in insurance networks. 
	\end{abstract}
\maketitle

%


\section{Introduction}

In this paper, we consider the following family of chance constrained optimization
problems:
\begin{align}
\label{chance-constraint-opt}%
\begin{array}
[c]{ll}%
\mbox{minimize} & c^{\top}x\\
\mbox{subject to} & \mathrm{P}(\phi(x,L) > 0)\leq\delta,\\
& x\in\mathbb{R}^{d_{x}}.
\end{array}
\tag{$\rm{CCP}_\delta$}%
\end{align}
where $x\in\mathbb{R}^{d_{x}}$ is a $d_{x}$-dimensional decision vector and $L$
is a $d_{l}$-dimensional random vector in $\mathbb{R}^{d_{l}}$. The elements
of $L$ are often referred to as risk factors; the function $\phi
:\mathbb{R}^{d_{x}}\times\mathbb{R}^{d_{l}}\rightarrow\mathbb{R}$ is often
assumed to be convex in $x$ and often models a cost constraint; the parameter
$\delta> 0$ is the risk level of the tolerance. \begin{revisionenv}
Our framework encompasses the joint chance constraint of the form $\mathrm{P}(\phi_i(x,L) > 0, ~ \exists i\in \{1,\ldots, n\} )\leq\delta$, by setting $\phi(x,L) = \max_{i = 1,\ldots, n}\phi_i(x,L)$.
\end{revisionenv}

Chance constrained optimization problems have a rich history in Operations
Research. Introduced by \cite{charnes1958cost}, chance constrained
optimization formulations have proved to be versatile in modeling and decision
making in a wide range of settings. For example,
\cite{prekopa1970probabilistic} used these types of formulations in the
context of production planning. The work of \cite{bonami2009exact} illustrates
how to take advantage of chance constrained optimization formulations in the
context of portfolio selection. In the context of power and energy control the
use of chance constrained optimization is illustrated in
\cite{andrieu2010model}. These are just examples of the wide range of
applications that have benefited (and continue to benefit) from chance
constrained optimization formulations and tools.

Consequently, there has been a significant amount
of research effort devoted to the solution of chance constrained optimization
problems. Unfortunately, however, these types of problems are provably NP-hard
in the worst case, see \cite{luedtke2010integer}. As a consequence, much of
the methodological effort has been placed into developing: a) solutions in the
case of specific models; b) convex and, more generally, tractable
relaxations;
c) combinatorial optimization tools; d) Monte-Carlo sampling schemes. \begin{revisionenv}
Of course, hybrid approaches are also developed. For example, as a combination of type b) and type d) approaches, \cite{hong2011sequential} show that the solution to a chance constraint optimization problem can be approximated by  optimization problems with constraints represented as the difference of two convex functions. In turn, this is further approximated by solving a sequence of convex optimization problems, each of which can be solved by a gradient based Monte Carlo method. Another example is \cite{pena2020solving}, which combines relaxations of type b) with sample-average approximation associated with type d) methods. In addition to the aforementioned types, \cite{hong2020learning} provides an upper bound for the chance constraint optimization problem using a robust optimization with a data-driven uncertainty set, achieving a dimension independent sample complexity.
\end{revisionenv}

Examples of type a) approaches include the study of Gaussian or elliptical
distributions when $\phi$ is affine both in $L$ and $x$. In this case, the
problem admits a conic programming formulation, which can be efficiently
solved, see \cite{lagoa2005probabilistically}.
\begin{revisionenv}
Type b) approaches include \cite{hillier1967chance,seppala1971constructing, ben2000robust, ben2002robust,prekopa2003probabilistic,bertsimas2004price,nemirovski2006convex, chen2010cvar, tong2020optimization}.
These approaches usually integrate probabilistic inequalities such as Chebyshev's bound, Bonferroni’s bound, Bernstein's approximations, or large deviation principles to 
construct tractable analytical approximations.
Type c) methods are based on branch and bounding algorithms, which connect
squarely with the class of tools studied in areas such as integer programming,
see \cite{ahmed2008solving, luedtke2010integer,
kuccukyavuz2012mixing,
luedtke2014branch, zhang2014branch, lejeune2016solving}.
Type d) methods include the sample gradient method, the sample average approximation and the scenario approach. The sample gradient method is usually combined with a smooth approximation, see \cite{hong2011sequential} for example. The sample average approximations studied by \cite{luedtke2008sample} and \cite{barrera2016chance}, although simplifying the constraint's probabilistic structure via replacing the population distribution by sampled empirical distribution, are nevertheless hard to solve due to non-convex feasible regions. \end{revisionenv}
The method we consider in this paper is the scenario approach. The scenario approach
is introduced and studied in \cite{calafiore2005uncertain} and is further
developed in a series of papers, including
\cite{calafiore2006scenario,nemirovski2006scenario}. 

The scenario approach is the most popular generic method for 
(approximately) solving chance constrained optimization. The idea is to sample a
number $N$ of scenarios \begin{revisionenv}
(each scenario consists of a sample of $L$) 
\end{revisionenv}
and enforce the constraint in all of these scenarios.
The intuition is that if for any scenario, say $L^{(i)}$, the constraint
$\phi(L^{(i)},x)<0$ is convex in $x$, and $\delta>0$ is small, we expect that by
suitably choosing $N$ the constrained regions can be relaxed by enforcing
$\phi(L^{(i)},x)<0$ for all $i = 1,\ldots, N$, leading to a good and, in some
sense, tractable (if $N$ is of moderate size) approximation of the chance
constrained region. Of course, this intuition is correct only when $\delta>0$
is small and we expect the choice of $N$ to be largely influenced by this
asymptotic regime.

By choosing $N$ sufficiently large, the scenario approach allows obtaining both upper and lower bounds which become asymptotically tighter as
$\delta\rightarrow0$. In a celebrated paper, \cite{calafiore2006scenario}
provide rigorous support for this claim. In particular, given a confidence
level $\beta\in(0,1)$, if $N\geq (2/\delta) \times\log(1/\beta) + 2d +
(2d/\delta)\times\log(2/\delta)$, with probability at least $1-\beta$, the
optimal solution of the scenario approach relaxation is feasible for the original chance constrained problem and, therefore, an upper bound to the
problem is obtained. 

Unfortunately, the required sample size of $N$ grows with $(1/\delta)\times\log(1/\delta)$ as $\delta$ becomes small, limiting the scope of the scenario methods in applications. 
\begin{revisionenv}
Many applications of chance constraint optimization require a very small $\delta$. For example, in the 5G ultra-reliable communication system design, the failure probability $\delta$ is no larger than $10^{-5}$, see \cite{alsenwi2019chance}; for fixed income portfolio optimization, an investment grade portfolio has a historical default rate of $10^{-4}$, reported by  \cite{frank2008municipal}.
\end{revisionenv}

Motivated by this, \cite{nemirovski2006scenario} developed a method that lowers the required
sample size to the order of $\log(1/\delta)$, making additional assumptions on
the function $\phi$ (which is taken to be bi-affine), and the risk factors
$L$, which are to be assumed light-tailed. Specifically, the moment generating
function $E[\exp(s L) ]$ is assumed to be finite in a neighborhood of the origin. No guarantee is given in terms of how far the upper bound is from the optimal value function of the problem as $\delta \rightarrow 0$.

In the present paper, we focus on improving the scalability of $N$ in terms of
$1/\delta$ for the practically important case of heavy-tailed risk factors.
Heavy-tailed distributions appear in a wide range of applications in science,
engineering and business, see e.g., \cite{embrechts2013modelling},
\cite{wierman2012tail}, but, in some aspects, are not as well understood as light-tails. One reason is that techniques from convex duality cannot be applied as the moment
generating function of $L$ does not exist in a neighborhood of $0$. In
addition, probabilistic inequalities, exploited in
\cite{nemirovski2006scenario}, do not hold in this setting. Only very
recently, a versatile algorithm for heavy-tailed rare event simulation has
been developed in \cite{chen2019efficient}.

The main contribution of our paper is an algorithm that provides a sample
complexity for $N$ which is bounded in $1/\delta$, assuming a versatile class
of heavy-tailed distributions for $L$. Specifically, we shall assume that $L$
follows a semi-parametric class of models known as multivariate regular
variation, which is quite standard in multivariate heavy-tail modeling, cf.\ 
\cite{embrechts2013modelling, resnick2013extreme}. A precise definition is
given in Section \ref{sec-construction}. Moreover, our estimator is shown to be within a constant factor to the solution to \eqref{chance-constraint-opt} with high probability, uniformly as $\delta \rightarrow 0$. 
\begin{revisionenv}
We are not aware of other approaches that provide a uniform performance guarantee of this type.
\end{revisionenv}


We illustrate our assumptions and our framework with a risk problem of
independent interest. This problem consists in computing a collective salvage
fund in a network of financial entities whose liabilities and payments are
settled in an optimal way using the Eisenberg-Noe model, see \cite{eisenberg2001systemic}. The salvage fund is computed to minimize its
size in order to guarantee a probability of collective default after
settlements of less than a small prescribed margin. 
\begin{revisionenv}
For the sake of demonstrating the broad applicability of our method, we present a portfolio optimization problem with value-at-risk constraints as an additional running example.
\end{revisionenv}

The rest of the paper is organized as follows.
In Section \ref{sec-running-examples}, we introduce the portfolio optimization problem and the minimal salvage fund problem as particular applications of chance constraint optimization. We employ both problems as running examples to provide a concrete and intuitive explanation for the concepts we introduce throughout the paper.
\begin{revisionenv}
In Section \ref{sec-review-scenario-approach}, we provide a brief review of the scenario approach in \cite{calafiore2006scenario}. 

The ideas behind our main algorithmic contributions are given in Section \ref{sec-algorithm}, where we introduce its intuition, rooted in ideas originating from rare event simulation. \end{revisionenv} Our algorithm requires the
construction of several auxiliary functions and sets. How to do this is
detailed in Section \ref{sec-construction}, in which we also present several
additional technical assumptions required by our constructions. In Section \ref{sec-construction}, we also explain that our procedure results in an estimate which is within a constant factor of the optimal solution of the underlying chance constrained problem with high probability as $\delta \rightarrow 0$. In Section \ref{sec-example} we show that the assumptions
imposed are valid in our motivating example (as well as a second example with quadratic cost structure inside the probabilistic constraint). Numerical results for the examples are provided in Section \ref{sec-numerics}. Throughout our discussion in each section we present a series of results which summarize the main ideas of our constructions. To keep the discussion fluid, we  present the corresponding proofs in Appendix \ref{sec-proofs} unless otherwise indicated.

\textbf{Notations:} in the sequel, $\mathbb{R}_{+} = [0,+\infty)$ is the set
of non-negative real numbers, $\mathbb{R}_{++} = (0,+\infty)$ is the set of
positive real numbers, and $\overline{\mathbb{R}} = [-\infty, +\infty]$ is the
extended real line. A column vector with zeros is denoted by $\mathbf{0}$,
and a column vector with ones is denoted by $\mathbf{1}$. For any matrix
$Q$, the transpose of $Q$ is denoted by $Q^{\top}$; the Frobenius norm of $Q$ is denoted by $\|Q\|_F$. The identity matrix is
denoted by $I$. For two column vectors $x, y \in\mathbb{R}^{d}$, we say
$x\preceq y$ if and only if $y-x\in\mathbb{R}^{d}_{+}$. For $\alpha
\in\mathbb{R}$ and $x \in\mathbb{R}^{d}$, we use $\alpha\cdot x$ to denote the
scalar multiplication of $x$ with $\alpha$. For $\alpha\in\mathbb{R}$ and
$E\subseteq\mathbb{R}^{d}$, we define $\alpha\cdot E = \{\alpha\cdot x\mid
x\in E\}$. The optimal value of an optimization problem $(\rm{P})$ is denoted by $\Val(\rm{P})$. We also use Landau's notation. In particular, if $f(\cdot)$ and $g(\cdot)$ are non-negative real valued functions, we write $f(t)=O(g(t))$ if $f(t)\leq c_0 \times g(t))$ for some $c_0 \in (0,\infty)$ and $f(t)=\Omega(g(t))$ if $f(t)\geq g(t))/c_0$ for some $c_0 \in (0,\infty)$.

\section{Running Examples}
\label{sec-running-examples}
\subsection{Portfolio Optimization with VaR Constraint}
\begin{revisionenv}
We first introduce a portfolio optimization problem. Suppose that there are $d$ assets to invest. If we invest a dollar in the $i$-th asset, the investment has mean return $\mu_i$ and a non-negative random loss $L_i$. Let $x = (x_1, \cdots, x_d)$ represent the amount of dollar invested in different assets, and let $\mu = (\mu_1,\ldots, \mu_d)$ and $L = (L_1, \ldots, L_d)$. We assume that $L$ follows a multivariate heavy-tailed distribution, in a way made precise later on.  The portfolio manager's goal is to maximize the mean return of the portfolio, which is equal to $\mu^{\top}x$, with a portfolio risk constraint prescribed by a risk measure called value-at-risk (VaR). The VaR at level $1-\delta \in (0,1)$ for a random variable $X$ is defined as
\begin{align*}
    \mathrm{VaR}_{1-\delta}(X)
    = \min\{z\in \mathbb{R}: F_X(z)\geq 1-\delta\}.
\end{align*}

For a given number $\eta>0$, we formulate the following portfolio optimization problem.
\begin{align*}
\begin{array}
[c]{ll}%
\mbox{maximize} & \mu^{\top}
x\\
\mbox{subject to} &  \mathrm{VaR}_{1-\delta}(x^{\top} L) 
\leq \eta,\\
& x\in\mathbb{R}_{++}^{d}.\\
\end{array}
\end{align*}
Using the definition of VaR and the fact that the cumulative distribution function is right continuous, we conclude that 
$\mathrm{VaR}_{1-\delta}(x^{\top} L)\leq \eta$ is equivalent to $\mathrm{P} (x^{\top}L - \eta > 0) \leq \delta$. In order to facilitate the technical exposition, we apply the change of variable $x_i \mapsto 1/x_i$ to homogenize the constraint function, yielding the following equivalent chance constrained optimization problem in standard form:
\begin{align}
\label{eq-portfolio-opt}
\begin{array}
[c]{ll}%
\mbox{maximize} & \sum_{i=1}^{d} (\mu_i/x_i) \\
\mbox{subject to} & \mathrm{P} \big(\phi(x,L)> 0\big) \leq \delta,\\
& x\in\mathbb{R}_{++}^{d}.
\\
\end{array}
\end{align}
where $\phi(x,l) = \sum_{i=1}^{d} (l_i/x_i) - \eta$. Despite the nonlinear objective, \cite[Section 4.3]{calafiore2005uncertain} shows that it admits an epigraphic reformulation with a linear objective so that the standard scenario approach is applicable.

\end{revisionenv}

\subsection{Minimal Salvage Fund}
\label{sec-salvage-fund}

Suppose that there are $d$ entities or firms, which we can interpret as
(re)insurance  firms. Let $L = (L_{1},\ldots, L_{d}) \in \mathbb{R}_{+}^d$ denotes the
vector of incurred losses by each firm, where $L_{i}$ denotes the total
incurred loss that entity $i$ is responsible to pay. We assume that $L$
follows a multivariate heavy-tailed distribution as in the previous example. Let $Q = (Q_{i,j}: i, j \in\{1, \ldots, d\})$ be a deterministic matrix
where $Q_{i,j}$ denotes the amount of money received by entity $j$ when entity
$i$ pays one dollar. We assume that $Q_{i,j} \geq0$ and $\sum_{j = 1}^{d}Q_{i,j} < 1$. 
\begin{revisionenv}
Let $x = (x_{1},\ldots,x_{d})$ denote the total amount that the salvage fund allocated to each entity, and $y^{\ast} = (y_{1}^{\ast}
,\ldots,y_{d}^{\ast} )$ denote the amount of the final settlement. The amount of final settlement is determined by the following optimization problem:
\begin{align*}
y^{\ast} = y^{\ast}(x,L) =\arg\max\{\mathbf{1}^{\top}y\mid0\preceq y\preceq L,
\quad\left(  I-Q^{\top}\right)  y\preceq x\}.
\end{align*}
\end{revisionenv}
In words, the system maximizes the payments subject to the constraint that
nobody pays more than what they have (in the final settlement), and nobody pays more
than what they owe. Notice that $y^{\ast}= y^{\ast}(x,L)$ is also a random
variable (the randomness comes from $L$) satisfying $\mathbf{0}\preceq
y^{\ast}\preceq L$. 
\begin{revisionenv}
Suppose that entity $i$ bankrupts if the deficit
$L_{i} - y_{i}^{\ast}\geq m_i$, where $m\in \mathbb{R}^{d}_{+}$ is a given vector.
\end{revisionenv} 
We are interested in finding the minimal amount of salvage fund that ensures no bankruptcy happens with probability at least $1-\delta.$ The problem can be
formulated as a chance constraint programming problem as follows
\begin{align}%
\begin{array}
[c]{ll}%
\mbox{minimize} & \mathbf{1}^{\top}x\\
\mbox{subject to} & \mathrm{P}(L - y^{\ast}(x,L) \preceq m)\geq
1-\delta,\\
& x\in\mathbb{R}_{++}^{d}.
\end{array}
\label{minimal-salvage-fund-problem}%
\end{align}
Now we write the problem \eqref{minimal-salvage-fund-problem} into standard
form. Notice that $L - y^{\ast}(x,L) \preceq m$ if and only if $\phi(x,L)
\leq 0$, where $\phi(x,L)$ is defined as follows
\begin{align*}
\phi\left( x,L\right)  &:=\min_{b,y}
\{b \mid(L - y - m) \preceq b\cdot\mathbf{1},\;
\left(  I-Q^{\top}\right)  y\preceq x, \; y\succeq\mathbf{0}\}.
\end{align*}
Therefore, problem \eqref{minimal-salvage-fund-problem} is equivalent to
\begin{align}%
\begin{array}
[c]{ll}%
\mbox{minimize} & \mathbf{1}^{\top}x\\
\mbox{subject to} & \mathrm{P}(\phi(x,L) > 0)\leq\delta,\\
& x\in\mathbb{R}_{++}^{d}.
\end{array}
\label{chance-constraint-optimization-general}%
\end{align}

\section{Review of Scenario Approach}
\label{sec-review-scenario-approach}
As mentioned in the introduction, a popular approach to solve the chance
constraint problem proceeds by using the scenario approach developed by
\cite{calafiore2006scenario}. They suggest to approximate the probabilistic
constraint $\mathrm{P}(\phi(x,L) > 0)\leq\delta$ by $N$ sampled constraints
$\phi(x,L^{(i)}) \leq0$ for $i = 1,\ldots,N$, where $\{L^{(1)}, \ldots
,L^{(N)}\}$ are independent samples. Instead of solving the original chance
constraint problem \eqref{chance-constraint-opt}, which is usually
intractable, we turn to solve the following optimization problem
\begin{align}
\label{sampled-problem}%
\begin{array}
[c]{ll}%
\mbox{minimize} & c^{\top}x\\
\mbox{subject to} & \phi(x,L^{(i)}) \leq0, \quad i = 1,\ldots,N,\\
& x\in\mathbb{R}^{d_{x}}.
\end{array}
\tag{$\rm{SP}_N$}%
\end{align}
The total sample size $N$ should be large enough to ensure the feasible
solution to the sampled problem \eqref{sampled-problem} is also a feasible
solution to the original problem \eqref{chance-constraint-opt} with a high
confidence level. According to \cite{calafiore2006scenario}, for any given
confidence level parameter $\beta\in(0,1)$, if
\begin{align*}
N\geq\frac{2}{\delta} \log\frac{1}{\beta} + 2d + \frac{2d}{\delta} \log\frac
{2}{\delta},
\end{align*}
then any feasible solution to the sampled optimization problem
\eqref{sampled-problem} is also a feasible solution to
\eqref{chance-constraint-opt} with probability at least $1-\beta$. However,
when $\delta$ is small, the total number of sampled constraints is of order
$\Omega((1/\delta)\log(1/\delta))$, which could be a problem for
implementation. For example, as we shall see in Section \ref{sec-numerics},
when $\beta= 10^{-5}$, $d = 15$ and $\delta= 10^{-3}$, the number of sampled
constraints $N$ is required to be larger than $2\times10^{5}$. 
In contrast, our method only requires to sample $2\times 10^3$ constraints. 
 
\section{General Algorithmic Idea}

\label{sec-algorithm}

To facilitate the development of our algorithm, we introduce some additional notation and a desired technical property. As we shall see, if the technical property is
satisfied, then there is a natural way to construct a scenario approach based
algorithm that only requires $O(1)$ of total sampled constraints. \begin{revisionenv} We exploit key intuition borrowed from rare event simulation. A common technique exploited, for example, in \cite{chen2019efficient}, is the construction of a so-called super set, which contains the rare event of interest. The super set should be constructed with a probability which is of the same order as that of the rare event of interest. If the conditional distribution given being in the super set is accessible, this can be used as an efficient sampling scheme. The first part of this section simply articulates the elements involved in setting the stage for constructing such a set in the outcome space of $L$. Later, in Section \ref{sec-construction}, we will impose assumptions in order to ensure that the probability of the super set, which eventually we will denote by $C_{\delta}$ is suitably controlled as $\delta \rightarrow 0$. Simply collecting the elements necessary to construct $C_{\delta}$ requires introducing some super sets involving the decision space, since the optimal decision is unknown. \end{revisionenv}

Let $F_{\delta}\subseteq\mathbb{R}^{d_{x}}$ denote the feasible region of the chance constraint optimization problem \eqref{chance-constraint-opt},
i.e., 
\begin{align}
F_{\delta} := \{x\in\mathbb{R}^{d_{x}}\mid\mathrm{P}(\phi(x, L) > 0)\leq
\delta\}.\label{eq-feasible-region}%
\end{align}
Here, the subscript $\delta$ is involved to emphasize that the feasible region
$F_{\delta}$ is parametrized by the risk level $\delta$. For any fixed
$x\in\mathbb{R}^{d_{x}}$, let $V_{x} := \{L\in\mathbb{R}^{d_{l}}\mid
\phi(x,L)>0\}$ denote the \emph{violation event at $x$}.

\begin{property}
	\label{property}  For any $\delta>0$, there exist a set $O_{\delta}%
	\subseteq\mathbb{R}^{d_{x}}$, and an event $C_{\delta}\subseteq\mathbb{R}%
	^{d_{l}}$ that satisfy the following statements. 
	
	\begin{enumerate}[label = \textrm{\alph*)}] 
		
		\item The feasible set $F_{\delta}$ is a subset of 
		$O_{\delta}$. 
		
		\item The event $C_{\delta}$ contains the violation event $V_{x}$ for any
		$x\in O_{\delta}$. 
		
		\item There exist a constant $M > 0$ independent of $\delta$ such that
		$\mathrm{P}(L\in C_{\delta}) \leq M\cdot\delta$. 
	\end{enumerate}
\end{property}

\begin{figure}[ptb]
	\centering
	\includegraphics[width = 5in]{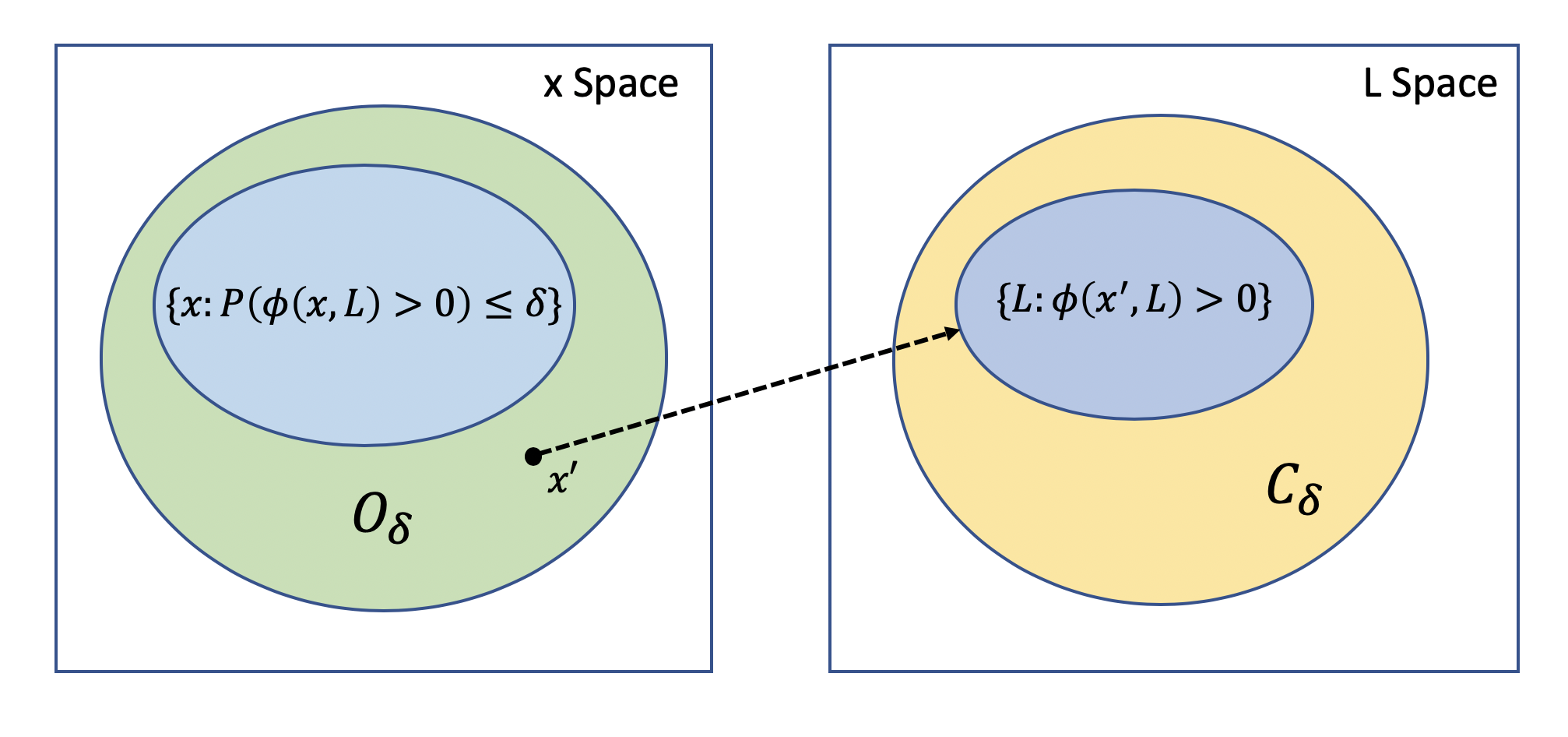}  \vspace{-0.3in}
	\caption{Pictorial illustration of $O_{\delta}$ and $C_{\delta}$.}%
	\label{fig:set_inclusion}%
\end{figure}

In the rest of this paper, we will refer to $O_{\delta}$ as \emph{the outer
	approximation set}, and $C_{\delta}$ as \emph{the uniform conditional event}.
A graphical illustration of $O_{\delta}$ and $C_{\delta}$ is shown in Figure
\ref{fig:set_inclusion}.

Now, given $O_{\delta}$ and $C_{\delta}$ that satisfies Property
\ref{property}, we define the conditional sampled
problem \eqref{conditional-sampled-problem}:
\begin{align}
\label{conditional-sampled-problem}
\begin{array}
[c]{ll}%
\mbox{minimize} & c^{\top}x\\
\mbox{subject to} & \phi(x, L^{(i)}_{\delta}) \leq0, \quad i = 1,\ldots
,N^{\prime}.\\
& x\in O_{\delta}.
\end{array}
\tag{$\rm{CSP}_{\delta,N'}$}%
\end{align}
where $L^{(i)}_{\delta}$ are i.i.d. samples generated from the conditional
distribution $(L|L\in C_{\delta})$.

\begin{revisionenv}
We now present our main result of this section in Lemma \ref{lemma-main}, which validates \eqref{conditional-sampled-problem} is an effective and sample efficient scenario approximation by incorporating \cite[Theorem 2]{calafiore2006scenario} and Property \ref{property}. The proof of Lemma \ref{lemma-main} will be presented in Section \ref{sec-proof-main}.
\end{revisionenv}

\begin{lemma} 
	\label{lemma-main}  Suppose that Property \ref{property} is imposed, and let
	$\beta> 0$ be a given confidence level. 
	\begin{enumerate}
		\item Let $\delta^{\prime} =\delta/\mathrm{P}(L\in C_{\delta})\geq1/M$ and $N'$ be any integer that satisfies
	    \begin{align}
    	\label{eq-bound-N-prime}N^{\prime}\geq
    	\frac{2}{\delta^{\prime}} \log\frac{1}{\beta} + 2d + \frac{2d}%
    	{\delta^{\prime}} \log\frac{2}{\delta^{\prime}}.
    	\end{align}
		With probability at least $1-\beta$, if the conditional sampled problem
		\eqref{conditional-sampled-problem} is feasible, then its optimal solution
		$x^{\ast}_{N}\in F_{\delta}$ and 
		$\Val\eqref{conditional-sampled-problem}\geq \Val \eqref{chance-constraint-opt}.$ 
		
		\item Let $N'$ be any integer such that $N'\leq \beta\delta^{-1}\mathrm{P}(L\in C_{\delta})$. Assume that the chance constraint
		problem ($CCP_{\delta}$) is feasible. Then, with probability at least
		$1-\beta$, $\Val \eqref{chance-constraint-opt}\geq\Val\eqref{conditional-sampled-problem}$. 
	\end{enumerate}
\end{lemma}

\begin{remark}
	Note that the lower bound given in \eqref{eq-bound-N-prime} is not greater than $2M \log(\frac{1}{\beta}) + 2d + 2dM\log(2M)$, which is independent of $\delta$. Therefore, Lemma \ref{lemma-main} shows that the chance constraint
	problem \eqref{chance-constraint-opt} can be approximated by \eqref{conditional-sampled-problem} with sample complexity bounded uniformly as  $\delta \rightarrow 0$, as long as Property \ref{property} is satisfied.
\end{remark}
\begin{remark}
     Efficiently generating samples of $(L|L\in C_\delta)$ when $\delta\rightarrow 0$ requires rare event simulation techniques.
     For example, when $L$ is light-tailed, exponential tilting can be applied to achieve $O(1)$ sample complexity uniformly in $\delta$; when $L$ is heavy-tailed, with the help of specific problem structure, one can apply importance sampling, see \cite{blanchet2010efficient}, or Markov Chain Monte Carlo, see \cite{gudmundsson2014markov}, to design an efficient sampling scheme. The specific structure of our salvage fund example results in $C_\delta$ being the complement of a box, which makes the sampling very tractable if the element of $L$ are independent. 
     
     Even if the aforementioned rare event simulation techniques are hard to apply in practice, we can still apply a simple acceptance-rejection procedure to sample the conditional distribution $(L|L\in C_\delta)$. It costs $O(1/\delta)$ samples of $L$ on average to get one sample of $(L|L\in C_\delta)$, since $\mathrm{P}(L\in C_{\delta}) = O(\delta)$. Consequently, the total complexity for generating $L^{(i)}_{\delta}, i = 1,\ldots, N'$ and solving $\eqref{conditional-sampled-problem}$ is $O(1/\delta)$, which is still much more efficient than the scenario approach in \cite{calafiore2006scenario}, because it requires computational complexity $O(((1/\delta)\log(1/\delta))^3)$ for solving a linear programming problem with $O((1/\delta)\log(1/\delta))$ sampled constraints by the interior point method.
\end{remark}

Although Property \ref{property} seems to be restrictive at  first
glance, we are still able to construct the sets $O_{\delta}$ and $C_{\delta}$
for a rich class of functions $\phi(x,L)$, including the constraint
function for the minimal salvage fund problem. As we shall see in the proof of
Lemma \ref{lemma-main}, once $O_{\delta}$ and $C_{\delta}$ are constructed the
sampled problem \eqref{sampled-problem-eq} is a tractable approximation to the
problem \eqref{chance-constraint-opt}. We  explain how to construct the
sets $O_{\delta}$ and $C_{\delta}$ in the next section under some additional
assumptions. These assumptions relate in particular to the distribution of
$L$. It turns out that, if $L$ is heavy-tailed, the construction of
$O_{\delta}$ and $C_{\delta}$ becomes tractable.
\subsection{Proof of Lemma \ref{lemma-main}.}\label{sec-proof-main}
	If Property \ref{property} is satisfied,  
	\eqref{chance-constraint-opt} is equivalent
	to
	\begin{align}
	\label{chance-constraint-opt-equi}%
	\begin{array}
	[c]{ll}%
	\mbox{minimize} & c^{\top}x\\
	\mbox{subject to} & \mathrm{P}(\phi(x, L) > 0 \mid L\in C_{\delta})\leq
	\delta/\mathrm{P}(L\in C_{\delta}),\\
	& x\in O_{\delta}\subseteq\mathbb{R}^{d_{x}}.
	\end{array}
	\end{align}
	Let $\delta^{\prime}:=\delta/\mathrm{P}(L\in C_{\delta})\geq1/M$ denote the risk
	level in the equivalent problem \eqref{chance-constraint-opt-equi}. The sampled optimization problem related to problem
	\eqref{chance-constraint-opt-equi} is given by
	\begin{align}
	\label{sampled-problem-eq}%
	\begin{array}
	[c]{ll}%
	\mbox{minimize} & c^{\top}x\\
	\mbox{subject to} & \phi(x, L^{(i)}_{\delta}) \leq0, \quad i = 1,\ldots
	,N^{\prime},\\
	& x\in O_{\delta},
	\end{array}
	\tag{\ref{conditional-sampled-problem}}
	\end{align}
	where the $L^{(i)}_{\delta}$ are independently sampled from $\mathrm{P}(\cdot\mid L\in
	C_{\delta})$. Notice that 
	\begin{align*}
    	N^{\prime}
    	\geq\frac{2}{\delta^{\prime}} \log\frac{1}{\beta} + 2d + \frac{2d}%
    	{\delta^{\prime}} \log\frac{2}{\delta^{\prime}}.
    \end{align*} 
    According to \cite[Corollary 1 and Theorem 2]{calafiore2006scenario}, with probability at least $1-\beta$, if the sampled problem \eqref{sampled-problem-eq} is feasible, then the optimal solution to problem
	\eqref{sampled-problem-eq} is feasible to the chance constraint problem \eqref{chance-constraint-opt-equi}, thus it is also feasible for \eqref{chance-constraint-opt}. The proof of the first statement is complete.
	
	Now we turn to prove the second statement. Note that the equivalence between \eqref{chance-constraint-opt} and \eqref{chance-constraint-opt-equi} is still valid, so it is sufficient to compare the optimal values of \eqref{chance-constraint-opt-equi} and \eqref{conditional-sampled-problem}. By applying \cite[Theorem 2]{calafiore2006scenario} again, we have with probability at least $1-\beta$ the value of \eqref{conditional-sampled-problem} is smaller or equal than the optimal value of 
	\begin{align}
	\label{chance-constraint-opt-equi2}%
	\begin{array}
	[c]{ll}%
	\mbox{minimize} & c^{\top}x\\
	\mbox{subject to} & \mathrm{P}(\phi(x, L) > 0 \mid L\in C_{\delta})\leq
	1-(1-\beta)^{1/N'},\\
	& x\in O_{\delta}\subseteq\mathbb{R}^{d_{x}}.
	\end{array}
	\end{align}
	The proof is complete by using $1-(1-\beta)^{1/N'}\geq \beta/N'\geq \frac{\delta}{\mathrm{P}(L\in C_\delta)}$. So, using $\Val$ for ``value of'', $\Val \eqref{chance-constraint-opt-equi2}\leq\Val \eqref{chance-constraint-opt-equi}=\Val \eqref{chance-constraint-opt}$.
\section{Constructing outer approximations and summary of the algorithm}

\label{sec-construction}

\begin{revisionenv} In this section, we come full circle with the intuition borrowed from rare event simulation explained at the beginning of Section \ref{sec-algorithm}. The scale-free properties of heavy-tailed distributions (to be reviewed momentarily) coupled with natural (polynomial) growth conditions (like the linear loss) given by the structure of the optimization problem, provide the necessary ingredients to show that the set $C_\delta$ has a probability which is of order $O(\delta)$.
\end{revisionenv} In this section, we present two methods for the construction of $O_{\delta}$ and
$C_{\delta}$ satisfying Property \ref{property}. We mostly focus on our ``scaling method" which is presented in Section
\ref{sec-scaling}, which is facilitated precisely by the scale-free property that we will impose on $L$. 
After showing the construction of the outer sets under the scaling method, we summarize the algorithm at the end of Section \ref{sec-scaling}. We supply a lower bound guaranteeing a constant approximation for the output of the algorithm in Section \ref{sec-Const-Approx}. Our second method for outer approximation constructions is summarized in Section
\ref{sec-jointly-convex}. This method is simpler to apply because is based on linear approximations, however, it is somewhat less powerful because it assume that $\phi(x,L)$ is jointly convex. 

\subsection{Scaling Method}

\label{sec-scaling}

We are now ready to state our assumption on the distribution of $L$. We assume
that the distribution of $L$ is of multivariate regular variation, a
definition that we review first. For background, we refer to \cite{resnick2013extreme}.
Let $\mathcal{M}_{+}(\overline{\mathbb{R}}^{d_{l}}\backslash\{\mathbf{0}\})$ denote all Radon measures on the space
$\overline{\mathbb{R}}^{d_{l}}\backslash\{\mathbf{0}\}$ (recall that a measure is Radon
if it assigns finite mass to all compact sets).
If $\mu_{n}(\cdot),\mu(\cdot)\in\mathcal{M}_{+}(\overline{\mathbb{R}}^{d_{l}%
}\backslash\{\mathbf{0}\})$, then $\mu_{n}$ converges to $\mu$ vaguely, denoted by
$\mu_{n}\vconverge\mu$, if for all compactly supported
continuous functions $f:\overline{\mathbb{R}}^{d_{l}}\backslash\{\mathbf{0}%
\}\rightarrow\mathbb{R}_{+}$,
\[
\lim_{n\rightarrow\infty}\int_{\overline{\mathbb{R}}^{d_{l}}\backslash\{\mathbf{0}%
	\}}f(x)\mu_{n}(dx)=\int_{\overline{\mathbb{R}}^{d_{l}}\backslash\{\mathbf{0}\}}%
f(x)\mu(dx).
\]
 $L$ is  \emph{multivariate regularly varying} with
\emph{limit measure} $\mu(\cdot)\in\mathcal{M}_{+}(\overline{\mathbb{R}%
}^{d_{l}}\backslash\{\mathbf{0}\})$ if
\[
\frac{\mathrm{P}(x^{-1}L\in\cdot)}{\mathrm{P}(\Vert L\Vert_{2}>x)}%
\vconverge\mu(\cdot),\qquad\text{as }x\rightarrow\infty.
\]

\begin{assumption}
	\label{assumption-regular-varying} $L$ is multivariate
	regularly varying with limit measure $\mu(\cdot)\in\mathcal{M}_{+}%
	(\overline{\mathbb{R}}^{d_{l}}\backslash\{\mathbf{0}\})$.
\end{assumption}

We give some intuition behind this definition. Write $L$ in terms of polar
coordinates, with $R$ the radius and $\Theta$ a random variable taking values on the unit sphere. The radius $R =\Vert L\Vert_{2}$ has a one-dimensional regularly varying tail (i.e. we can write $\mathrm{P} (R>x) = L(x) x^{-\alpha}$ for a slowly varying function $L$ and $\alpha>0$). The angle $\Theta$, conditioned on $R$ being large, converges weakly (as $R\rightarrow\infty$) to a limiting random variable. The distribution of this limit can be expressed in terms of the measure $\mu$. For another recent application of multivariate regular variation in operations research, see \cite{kley2016risk}.

We proceed to analyze the feasible region $F_{\delta}$ when $\delta
\rightarrow0$. Intuitively, if the violation probability $\mathrm{P}(\phi(x,L)
> 0)$ has a strictly positive lower bound in any compact set, then $F_{\delta
}$ will ultimately be disjoint with the compact set when $\delta\rightarrow0$.
Thus, the set $F_{\delta}$ is expelled to infinity when $\delta\rightarrow0$ in
this case. $F_{\delta}$ is moving towards the direction that $\phi(x,L)$
becomes small such that the violation probability becomes smaller. For
instance, if $x$ is one dimensional and $\phi(x,L)$ is increasing in $x$, then $F_{\delta}$ is moving towards the negative direction. Consider
the portfolio optimization problem as another example, in which $\min_{i=1}^{d} x_{i}\rightarrow+\infty$ as
$\delta\rightarrow0$.

Now we begin to construct the outer approximation set $O_{\delta}$. To this
end, we need to introduce an auxiliary function which we shall call a
\emph{level function}. 

\begin{definition}\label{def-level-set}
	We say that $\pi:\mathbb{R}^{d_{x}}%
	\rightarrow\lbrack0,+\infty]$ is a level function if 
	\begin{enumerate}
		\item for any $\alpha\geq0$ and $x\in\mathbb{R}^{d_{x}}$, we have $\pi
		(\alpha\cdot x)=\alpha\cdot\pi(x)$, 
		
		\item $\lim_{\delta\rightarrow0}\inf_{x\in F_{\delta}}\pi(x)=+\infty$. 
	\end{enumerate}
	
	We also define the \emph{level set} $\Pi=\{x\in\mathbb{R}^{d_{x}}\mid
	\pi(x)=1\}$.
\end{definition}

As $F_{\delta}$ is moving to infinity, the level function is helpful to
characterize the `moving direction' of $F_{\delta}$ as well as the correct
rate of scaling as $\delta$ becomes small. As we shall see in the proof of
Lemma \ref{lemma-outer-approx}, for any $\delta$ small enough we can choose
some $\alpha_{\delta}$ and define
\[
O_{\delta}:=\bigcup_{\alpha\geq\alpha_{\delta}}\left(  \alpha\cdot\Pi\right)
\supseteq F_{\delta}.
\]
To construct $O_{\delta}$, we first select the level set $\Pi$, and then derive the scaling rate of $\alpha_{\delta}$.

The level function $\pi$ and the shape of $\Pi$ should be chosen in
accordance with the moving direction of $F_{\delta}$ to reduce the size of
$O_{\delta}$, in order to achieve better sample complexity. For example, when
$\phi(x,L)=-\Vert x\Vert^{2}-L$, the level function $\pi$ can be chosen as the
Euclidean norm and $\Pi$ can be chosen as the unit sphere in $\mathbb{R}%
^{d_{x}}$. For the portfolio optimization problem, the level function can be
chosen as $\pi(x)=\min_{i=1}^{d}x_{i}+\infty\cdot I(x\notin\mathbb{R}%
_{++}^{d_{x}})$ in accordance with our intuition that $\min_{i=1}^{d}%
x_{i}\rightarrow\infty$, and the level set can be chosen as $\Pi=\{x\in\mathbb{R}^{d_{x}}\mid
\min_{i=1}^{d}x_{i}=1\}$. Therefore, it is natural to impose the following assumption about the existence of the level function.

\begin{assumption}\label{assumption-level-set}
	There exist a level function $\pi$ and a level set $\Pi$.
\end{assumption}

To analyze the asymptotic shape of the uniform conditional event
$C_{\delta}$, we  connect the asymptotic distribution of $L$ to the
asymptotic distribution of $\phi(x,L)$. We pick a continuous non-decreasing
function $h:\mathbb{R}_{++}\rightarrow\mathbb{R}_{++}$ such that $\lim
_{\alpha\rightarrow+\infty}h(\alpha)=+\infty$ to characterize the scaling rate
of $L$. In addition, we pick another positive function $r:\mathbb{R}%
_{++}\rightarrow\mathbb{R}_{++}$ to characterize the scaling rate of
$\phi(\alpha\cdot x,h(\alpha)\cdot L)$. Intuitively, the scaling function
$r(\cdot)$ and $h(\cdot)$ should ensure the condition that $\{\frac
{1}{r(\alpha)}\phi(\alpha\cdot x,h(\alpha)\cdot L)\}_{\alpha\geq1}$ is tight.
For the minimal salvage fund problem with fixed $\delta$, as the
deficit $\phi(x,L)$ is asymptotically linear 
with respect to the salvage fund
$x$ and the loss $L$, we can simply pick $r(\alpha)=h(\alpha)=\alpha$ in this
problem. We next introduce two auxiliary functions $\Psi_{+}$ and $\Psi_{-}$.

\begin{definition}
	\label{def-asymp-bound}  Let $\Psi_{+}: \mathbb{R}^{d_{l}}\rightarrow
	\mathbb{R}$, $\Psi_{-}: \mathbb{R}^{d_{l}}\rightarrow\mathbb{R}$ be two Borel
	measurable functions. We say $\Psi_{+}$ (resp. $\Psi_{-}$) is the asymptotic
	uniform upper (resp. lower) bound of $\frac{1}{r(\alpha)}\phi(\alpha\cdot x,
	h(\alpha)\cdot l)$ over the level set $x\in\Pi$ if for any compact set
	$K\subseteq\mathbb{R}^{d_{l}}$,
	\begin{subequations}
		\label{eq-asymp-upper}%
		\begin{equation}
		\liminf_{\alpha\rightarrow\infty} \inf_{l\in K}\left( \Psi_{+}(l) - \sup
		_{x\in\Pi}\left[ \frac{1}{r(\alpha)}\phi(\alpha\cdot x, h(\alpha)\cdot
		l)\right] \right)  \geq0,
		\end{equation}
		\begin{equation}
		\label{eq-asymp-lower}\limsup_{\alpha\rightarrow\infty} \sup_{l\in K}\left(
		\Psi_{-}(l) - \inf_{x\in\Pi}\left[ \frac{1}{r(\alpha)}\phi(\alpha\cdot x,
		h(\alpha)\cdot l)\right] \right)  \leq0.
		\end{equation}
		
	\end{subequations}
\end{definition}

In Section \ref{sec-example}, we show for the salvage fund example how $\Psi_{+}$ and $\Psi_{-}$ can be written as maxima or minima of affine functions.
Here, we employ the functions $\Psi_{+}$ and $\Psi_{-}$ to define the event
$C_{\varepsilon,-}$ and $C_{\varepsilon,+}$, which serve as the inner and
outer approximation of the event $\cup_{x\in\Pi} V_{x}$, where $V_{x} =
\{l\in\mathbb{R}^{d_{l}}\mid\phi(x,l)>0\}$ is the violation event at $x$.

\begin{definition}
	\label{def-approx-event}  For $\varepsilon>0$, let $C_{\varepsilon,+}$
	(resp. $C_{\varepsilon,-}$) be the $\varepsilon$-outer (resp. inner)
	approximation event 
	\begin{subequations}
		\begin{equation}
		C_{\varepsilon,+} := \{l\in\mathbb{R}^{d_{l}}\mid\Psi_{+}(l) \geq
		-\varepsilon\},
		\end{equation}
		\begin{equation}
		C_{\varepsilon,-} := \{l\in\mathbb{R}^{d_{l}}\mid\Psi_{-}(l) \geq
		+\varepsilon\}.
		\end{equation}
		
	\end{subequations}
\end{definition}

We now define  $O_{\delta}:=\bigcup_{\alpha
	\geq\alpha_{\delta}} \alpha\cdot\Pi$. The following property ensures that the
shape of $\Pi$ is appropriate and $\alpha_{\delta}$ is large enough, hence
$O_{\delta}$ is an outer approximation of $F_{\delta}$.

\begin{property}
	\label{condition-outer-approx}  There exist $\delta_{0}$ such that for any
	$\delta<\delta_{0}$, we have an explicitly computable constant $\alpha_{\delta}$
	that satisfies
	\begin{align*}
	\mathrm{P}(\|L\|_{2}>h(\alpha_{\delta})) = O(\delta) \qquad\mbox{and}\qquad
	F_{\delta}\subseteq\bigcup_{\alpha\geq\alpha_{\delta}}\alpha\cdot\Pi=
	O_{\delta}.
	\end{align*}
	
\end{property}

If the violation probability is easy to analyze, we will directly derive the
expression of $\alpha_{\delta}$ and verify Property
\ref{condition-outer-approx}. Otherwise, we resort to Lemma
\ref{lemma-outer-approx}, which provides a sufficient condition of Property
\ref{condition-outer-approx} by analyzing the asymptotic distribution. The proof of Lemma \ref{lemma-outer-approx} is deferred to Appendix \ref{sec-proofs}.

\begin{lemma}
	\label{lemma-outer-approx} Suppose that Assumptions
	\ref{assumption-regular-varying} and \ref{assumption-level-set} hold.
	If there exists an asymptotic uniform lower bound function $\Psi_{-}(\cdot)$ as
	given in \eqref{eq-asymp-lower} and $\varepsilon>0$ such that $\mu
	(C_{\varepsilon,-})>0$, then Property \ref{condition-outer-approx} is satisfied.
\end{lemma}

We impose the following Assumption \ref{assumption-approx-func} on the
asymptotic uniform upper bound $\Psi_{+}(\cdot)$ so that we can employ the
multivariate regular variation of $L$  to estimate $\mathrm{P}(L \in
\alpha\cdot C_{\varepsilon, +})$ for large scaling factor $\alpha$.

\begin{assumption}
	\label{assumption-approx-func}  There exist an event $S\subseteq
	\mathbb{R}^{d_{l}}$ with $\mu(S^{c})<\infty$ such that
	\begin{align*}
	S\subseteq\alpha\cdot S, \qquad\Psi_{+}(l) \leq\Psi_{+}(\alpha\cdot l),
	\qquad\forall l\in S, \;\alpha\geq1.
	\end{align*}
	In addition, there exist some $\varepsilon>0$ such that $C_{\varepsilon, +}$
	is bounded away from the origin, i.e., $\inf_{l\in C_{\varepsilon, +}}\|l\|_{2}
	> 0.$
\end{assumption}

For the minimal salvage fund problem, since the deficit function $\phi(x,L)$
is coordinatewise nondecreasing with respect to the loss vector $L$, it is
reasonable to assume that its asymptotic bound $\Psi_{+}(\cdot)$ is also
coordinatewise nondecreasing. For this example, the closed form expression of
$\Psi_{+}(\cdot)$ and the detailed verification of all the assumptions are
deferred to Proposition \ref{prop-minimal-salvage-fund}. Our next result summarizes the construction of the outer approximation sets. 

\begin{theorem}
	\label{thm-set-construction}  Suppose that Property
	\ref{condition-outer-approx} and Assumption \ref{assumption-approx-func} are
	imposed. Then there exist $\delta_{0} > 0$ such that the following sets
	\begin{align}
	\label{eq-O-C-def}O_{\delta} = \bigcup_{\alpha\geq\alpha_{\delta}}\alpha
	\cdot\Pi, \qquad C_{\delta} = h(\alpha_{\delta})\cdot\big(C_{\varepsilon, +}
	\cup K^{c}\cup S^{c}\big)
	\end{align}
	satisfy Property \ref{property} for all $\delta<\delta_{0}$. Here, $S$ is
	given in Assumption \ref{assumption-approx-func} and $K$ is a ball in
	$\mathbb{R}^{d_{l}}$ with $\mu(K^{c}) <\infty$.
\end{theorem}

With the aid of Lemma \ref{lemma-main} and Theorem \ref{thm-set-construction}, we provide Algorithm \ref{algo-scenario} for approximating \eqref{chance-constraint-opt} in which the sampled optimization problem is bounded in $1/\delta$.

\begin{algorithm}[H]
\SetAlgoLined
 \SetKwInOut{Input}{input}
 \SetKwInOut{Output}{output}
 \Input{Constraint function $\phi$, risk tolerance parameter $\delta$, confidence level $\beta$, all the elements and  constants appearing in Property
	\ref{condition-outer-approx} and Assumption \ref{assumption-approx-func}.}
 Compute the expression of sets $O_\delta$ and $C_\delta$ by \eqref{eq-O-C-def}\;
 
 Compute required number of samples $N'$ by \eqref{eq-bound-N-prime}\;
 
 \For{$i = 1,\ldots, N'$}{
 Sample $L^{(i)}_{\delta}$ using acceptance-rejection or importance sampling.
 } 
 
 Solve the conditional sampled problem \eqref{conditional-sampled-problem}.
 
 \caption{Scenario Approach with Optimal Scenario Generation}
 \label{algo-scenario}
\end{algorithm}

In Section \ref{sec-Const-Approx}, our objective is to show that the output of the previous algorithm is guaranteed to be within a constant factor of the optimal solution to \eqref{chance-constraint-opt} with high probability, uniformly in $\delta$.

\subsection{Constant Approximation Guarantee}

\label{sec-Const-Approx}
We shall work under the setting of Theorem \ref{thm-set-construction}, so we enforce Property \ref{condition-outer-approx} and Assumptions \ref{assumption-approx-func}. We want to show that there exist some constant $\Lambda>1$ independent of $\delta$, such that $\Val \eqref{chance-constraint-opt}\leq \Val \eqref{conditional-sampled-problem} \leq \Lambda \times \Val \eqref{chance-constraint-opt}$ with high probability. This indicates that our result guarantees a constant approximation to \eqref{chance-constraint-opt} for regularly varying distributions (under our assumptions) in $O(1)$ sample complexity when $\delta \rightarrow 0$ with high probability.


Note that $\eqref{conditional-sampled-problem} \leq \Lambda \times \Val \eqref{chance-constraint-opt}$ is meaningful only if $\Val \eqref{chance-constraint-opt}>0$. We assume that the outer approximation set is good enough such that the following natural assumption is valid.
\begin{assumption}
    There exist $\delta>0$ such that $\min_{x\in O_\delta}c^{\top} x >0$.
    \label{assumption-positive-opt-val}
\end{assumption}

The previous assumption will typically hold if $c$ has strictly positive entries. 
Theorem \ref{thm-set-construction} and the form of $O_{\delta}$ guarantee that the norm of the optimal solution of \eqref{conditional-sampled-problem} grows in proportion to $\alpha_\delta$, so we also assume the following scaling property for $\phi(x,l)$. 
\begin{assumption}
    There exist a function $\phi_{\lim}:(\mathbb{R}^{d_{x}}\backslash\{\mathbf{0}\})\times (\mathbb{R}^{d_{l}}\backslash\{\mathbf{0}\})\rightarrow
    \mathbb{R}$ such that for every compact set $E\subseteq\mathbb{R}^{d_{l}}\backslash\{\mathbf{0}\}$, we have 
    \begin{align*}
    \lim_{\alpha\rightarrow \infty}
    \sup_{l\in E}
    \left\vert
    \frac{1}{r(\alpha)}\phi(\alpha\cdot x, h(\alpha)\cdot l)
    - \phi_{\lim}(x,l)
    \right\vert = 0.
    \end{align*}
    In addition, $\phi_{\lim}(x,l)$ is continuous in $l$. 
    \label{assumption-lim-phi}
\end{assumption}

\begin{revisionenv}
Assumption \ref{assumption-lim-phi} is satisfied by both running examples. For the portfolio optimization problem, we have 
$\phi(x,l) = \sum_{i=1}^{d} (L_i/x_i) - \eta$, thus $\phi_{\lim}(x,l) = \phi(x,l)$. For the minimal salvage fund problem, we have $\phi_{\lim}(x,l) = \phi(x,l) - m$ such that $\left\vert\alpha^{-1}\phi(\alpha\cdot x, \alpha\cdot l)
- \phi_{\lim}(x,l)\right\vert\leq \alpha^{-1}m$ and 
$|\phi_{\lim}(x,l)-\phi_{\lim}(x,l')|\leq \|l-l'\|_1$. 

\end{revisionenv}

We define the following optimization problem, which will serve as an asymptotic upper bound of $\eqref{conditional-sampled-problem}$ in stochastic order when $\delta\rightarrow 0$:
\begin{align}
	\label{limit-problem}
	\begin{array}
	[c]{ll}%
	\mbox{minimize} & c^{\top}x\\
	\mbox{subject to} & \phi_{\lim}(x, L_{\lim}^{(i)}) \leq0, \quad i = 1,\ldots
	,N^{\prime},\\
	& x\in \bigcup_{\alpha\geq1}\alpha
	\cdot\Pi,
	\end{array}
	\tag{$\rm{CSP}_{\lim,N'}$}
\end{align}
where $L_{\lim}^{(i)}$ are i.i.d. samples from a random variable $L_{\lim}$, whose distribution is characterized by $\mathrm{P}(L_{\lim}\in (C_{\varepsilon, +}
\cup K^{c}\cup S^{c})) = 1$ and $\mathrm{P}(L_{\lim}\in E) = \mu(E)/\mu(C_{\varepsilon, +}
\cup K^{c}\cup S^{c})$ for all measurable set $E \subseteq C_{\varepsilon, +}
\cup K^{c}\cup S^{c}$.

\begin{theorem}
	Let $\beta >0$ be a given confidence level and $N'$ be a fixed integer that satisfies \eqref{eq-bound-N-prime}. If Assumptions \ref{assumption-positive-opt-val} and \ref{assumption-lim-phi} are enforced, and \eqref{limit-problem} satisfies Slater's condition with probability one, then there exist $\delta_0 >0$ and $\Lambda>0$ such that
	\begin{equation*}
		\mathrm{P}\Big(
		\Val \eqref{chance-constraint-opt} \leq 
		\Val \eqref{conditional-sampled-problem}
		\leq \Lambda \times \Val \eqref{chance-constraint-opt}
		\Big) \geq 1-2\beta,
		\qquad\forall\delta<\delta_0.
	\end{equation*}
	\label{thm-nearly-optimal}
\end{theorem}
Slater's condition 
\begin{revisionenv}
(See Section 5.2.3 in \cite{boyd2004convex} for reference)
\end{revisionenv} 
can be verified directly on the problem \eqref{limit-problem}. This condition is satisfied in the salvage fund problem by standard linear programming duality. 

\subsection{Linear Approximation Method}
\begin{revisionenv}
\label{sec-jointly-convex} 
Suppose that the constraint function $\phi(x, l)$
is jointly convex in $(x,l)$, and $L$ is multivariate regularly varying.
We will develop a simpler method in this section to construct the outer
approximation set $O_{\delta}$ and the uniform conditional event $C_{\delta}$.

We first introduce a crucial assumption in the construction of $O_{\delta}$ and $C_{\delta}$.

\begin{assumption}
    \label{assumption-linear-approx}
    There exist a convex piecewise linear function $\phi_{-}(x, l):\mathbb{R}^{d_{x}}%
	\times\mathbb{R}^{d_{l}}\rightarrow\mathbb{R}$ of the form 
	\begin{align*}
		\phi_{-}(x, l) = \max_{i = 1,\ldots,N} a_{i}^{\top}l + b_{i}^{\top} x + c_{i}, \qquad
		a_{i}\in\mathbb{R}^{d_{l}}, b_{i}\in\mathbb{R}^{d_{x}}\mbox{ and } c_{i}\in\mathbb{R}\mbox{ for }i = 1,\ldots,N.
	\end{align*}
	such that:
	\begin{enumerate}
		\item $\phi_{-}(x,l) \leq\phi(x,l), ~\forall(x,l)\in\mathbb{R}^{d_{x}}%
		\times\mathbb{R}^{d_{l}}$;
		\item there exist some constant $C\in\mathbb{R}_{+}$ such that $\phi(x,l)\leq 0$ if $\phi_{-}(x, l)\leq -C$.
	\end{enumerate}
\end{assumption}

If $\phi(x,l)$ itself is a piecewise affine function, then Assumption \ref{assumption-linear-approx} is satisfied by simply taking $\phi_{-}(x,l) = \phi(x,l)$. For general jointly convex functions, the following lemma verifies Assumption \ref{assumption-linear-approx} if $\phi(x,l)$ has a compact zero sublevel set. 

\begin{lemma}
	\label{lemma-linear-approx}  If the constraint function $\phi(x, L): \mathbb{R}%
	^{d_{x}}\times\mathbb{R}^{d_{l}}\rightarrow\mathbb{R}$ is convex and twice
	continuously differentiable, and it has a compact zero
	sublevel set $Z_{\phi}:=\{(x,l)\in\mathbb{R}^{d_{x}}\times\mathbb{R}^{d_{l}}
	\mid\phi(x,l)\leq0\}$, then Assumption \ref{assumption-linear-approx} is satisfied.
\end{lemma}

With Assumption \ref{assumption-linear-approx} enforced, we are now ready to provide our main result in this section to fully summarize the construction of $O_{\delta}$ and
$C_{\delta}$.

\begin{theorem}
	\label{thm-convex}  If  Assumptions \ref{assumption-regular-varying}
	and \ref{assumption-linear-approx} hold, we can construct $O_{\delta}$
	and $C_{\delta}$ that satisfy Property \ref{property} as
	\begin{equation*}
	O_{\delta} := \bigcap_{i = 1}^{N} \{x\in\mathbb{R}^{d_{x}}\mid b_{i}^{\top}x +
	c_{i} + \bar{F}^{-1}_{a_{i}^{\top}L}(\delta)\leq0\}, \qquad C_{\delta}
	:=\bigcup_{i = 1}^{N} \{L\in\mathbb{R}^{d_{l}}\mid a_{i}^{\top}L + C > \bar
	{F}^{-1}_{a_{i}^{\top}L}(\delta)\},
	\end{equation*}
	where ${F}^{-1}_{a_{i}^{\top}L}(\delta)
	=\inf \{x \in \mathbb{R} \mid 
	\mathrm{P} (x > a_{i}^{\top}L)\leq \delta \}
	$.
\end{theorem}
\end{revisionenv}

\section{Verifying the Assumptions in Examples}
\label{sec-example}

In this section, we verify the elements required to apply our algorithm. We provide explicit expressions for sets $O_\delta$ and $C_\delta$ in the statement of the propositions. The detailed verification process and the steps for constructing sets $O_\delta$ and $C_\delta$ are presented as the proofs in 
Appendix \ref{sec-proofs}.

\subsection{Portfolio Optimization with VaR Constraint}\label{sec-verify-port-opt}

\begin{revisionenv}
In this section, we will verify that Theorem \ref{thm-set-construction} is applicable to an equivalent form of the portfolio optimization problem \eqref{eq-portfolio-opt}. 

\begin{proposition}
	\label{prop-port-opt}  The portfolio optimization problem
	\eqref{eq-portfolio-opt} satisfies all assumptions required by Theorem \ref{thm-set-construction}, such that the sets $O_\delta$ and $C_\delta$ admits the explicit expressions
	$$
	O_\delta = \left\{
    x\in\mathbb{R}_{++}^{d}~
    \middle\vert~ \eta\cdot x  \succeq \bar{F}^{-1}_{\mathbf{1}^{\top}L}(\delta)
    \right\},
    \qquad
    C_{\delta}
    = \left\{l\in\mathbb{R}_{++}^d ~\middle\vert~ 2\cdot \mathbf{1}^{\top} l \geq \bar{F}^{-1}_{\mathbf{1}^{\top}L}(\delta)
    \right\}.
	$$
\end{proposition}
\end{revisionenv}

\subsection{Minimal Salvage Fund}
\begin{revisionenv}
The key observation to solve the minimal salvage fund problem \eqref{chance-constraint-optimization-general} is the following lemma, which provides a closed form piecewise linear expression for the constraint function $\phi(x,L)$.
\begin{lemma}
    \label{lemma-minimal-salvage-fund-closed-form}
    In the minimal salvage fund problem \eqref{chance-constraint-optimization-general}, we have
    \begin{equation*}
        \phi(x, L) = \max_{i = 1,\ldots, d}\; L_i - \mathbf{e}^{\top}_{i}(I-Q^{\top})^{-1}x - m_i,
    \end{equation*}
    where $\mathbf{e}_{i}$ denote the unit vector on the $i$-th coordinate.
\end{lemma}

Now we prove that Theorem \ref{thm-convex} is applicable to the minimal salvage fund problem
\eqref{chance-constraint-optimization-general}.

\begin{proposition}
	\label{prop-minimal-salvage-fund} The minimal salvage fund problem
	\eqref{chance-constraint-optimization-general} satisfies all assumptions
	required by Theorem \ref{thm-convex}, such that the sets $O_\delta$ and $C_\delta$ admits the explicit expressions
	\begin{equation*}
	O_{\delta} = \bigcap_{i = 1}^{d} \{x\in\mathbb{R}^{d}\mid \bar{F}^{-1}_{L_i}(\delta)\leq
	\mathbf{e}^{\top}_{i}(I-Q^{\top})^{-1}x + m_i\}, \qquad C_{\delta}
	=\bigcup_{i = 1}^{d} \{L\in\mathbb{R}^{d}\mid L_i > \bar
	{F}^{-1}_{L_i}(\delta)\}.
	\end{equation*}
\end{proposition}
\end{revisionenv}

\subsection{Quadratic Model}

In this section, we consider a model with a quadratic control term in $x$ as an additional example.
Suppose that the constraint function $\phi(x,L): \mathbb{R}^{d_{x}}
\times\mathbb{R}^{d_{l}}\rightarrow\mathbb{R}$ is defined as
\begin{align}
\label{eq-constraint-quad}\phi(x, L) = x^{\top} Q x + x^{\top} A L,
\end{align}
where $Q\in\mathbb{R}^{d_{x}\times d_{x}}$ is a symmetric matrix and
$A\in\mathbb{R}^{d_{x}\times d_{l}}$ is a matrix with $\mathrm{rank}(A) =
d_{x}$, i.e. there exist $\sigma>0$ such that $\|A^{\top}x\|_2\geq\sigma\|x\|_2$. 

\begin{proposition}\label{prop-quadratic}
	Consider the chance constraint optimization model with constraint function
	defined as \eqref{eq-constraint-quad}. 
	
	\begin{enumerate}

		\item If $Q$ is a positive semi-definite matrix and $L$ has a positive
		density,  there exist some $\delta$ such that the problem is infeasible. 
		
		\item If $Q$ has a negative eigenvalue and $L$ is multivariate regularly
		varying,  the model satisfies all the assumptions required by Theorem
		\ref{thm-set-construction}. 
	\end{enumerate}
\end{proposition}

\section{Numerical Experiments}
\label{sec-numerics}

\begin{revisionenv}

In order to empirically study the computational complexity and compare the quality of the solutions, in this section we conduct numerical experiments for two scenario generation algorithms:
\begin{enumerate}
    \item the efficient scenario generation approach proposed in this paper (abbreviated as Eff-Sc);
    \item the scenario approach in \cite{calafiore2006scenario} (abbreviated as CC-Sc).
\end{enumerate}

In Section \ref{sec-numerics-portfolio-opt}, we present the results for the portfolio optimization problem. In Section \ref{sec-numerics-minimal-salvage-fund}, we present the results for the minimal salvage fund problem. The numerical experiment is conducted using a Laptop with 2.2 GHz Intel Core
i7 CPU, and the sampled linear programming problem is solved using CVXPY (\cite{cvxpy}) with the MOSEK solver (\cite{mosek}).  

\subsection{Portfolio Optimization with VaR Constraint}
\label{sec-numerics-portfolio-opt}
First of all, we present the parameter selection and the implement details for the numerical experiment of portfolio optimization problem \eqref{eq-portfolio-opt}. Suppose that there are $d = 10$ assets to invest, and the parameters of the problem are chosen as follows:
\begin{itemize}
    \item The mean return vector is $\mu = (1.0,1.5,2.0,2.5,3,1.6,1.2,1.1,1.8,2.2)$.
    \item $L_i$ are i.i.d. with Pareto cumulative distribution function $\mathrm{P}(L_i>l) = (l_i /l)$, for $l\geq\ell_i$.
    \item $\ell = (\ell_1,\ldots,\ell_d) = (2.1,1.3,1.6,2.5, 2.7, 1.3, 1.9, 1.5, 2.2, 2.3)$.
    \item The loss threshold $\eta = 1000$.
\end{itemize}

Now we explain the implementation detail of Eff-Sc. Recall the expression of $O_\delta$ and $C_\delta$ from Proposition \ref{prop-port-opt}, which involves the analytically unknown quantity $\bar{F}^{-1}_{\mathbf{1}^{\top}L}(\delta)$. Since quantile estimation is much more computationally efficient than solving the sampled optimization problem, we generate samples of $L$ to estimate a confidence interval of $\bar{F}^{-1}_{\mathbf{1}^{\top}L}(\delta)$ with large enough confidence level $1-o(\beta)$, and we denote the resulting confidence interval by $(\widehat{\text{LB}}, \widehat{\text{UB}})$. We replace the expressions of $O_\delta$ and $C_\delta$ by their sampled version conservative approximations, i.e.
$$
	O_\delta = \left\{
    x\in\mathbb{R}_{++}^{d}~
    \middle\vert~ \eta\cdot x  \succeq \widehat{\text{UB}}
    \right\},
    \qquad
    C_{\delta}
    = \left\{l\in\mathbb{R}_{++}^d ~\middle\vert~ 2\cdot \mathbf{1}^{\top} l \geq \widehat{\text{LB}}
    \right\}.
$$
The value of $\mathrm{P}(L\in C_\delta)$ is also estimated using the generated samples. We compute the required number of samples $N'$ using Lemma \ref{lemma-main}, and the samples of $L_{\delta}$ is generated via acceptance-rejection.

In Figure \ref{fig-efficiency-portfolio-opt}, we compare the efficiency between Eff-Sc and CC-Sc. Figure \ref{subfig-portfolio-opt-num-samples} presents the required number of samples for both algorithms, in which one can quickly remark that Eff-Sc requires significantly fewer samples than CC-Sc, especially for the problems with small $\delta$. In Figure \ref{subfig-portfolio-opt-cpu-time} we compare the running time for both models. Whereas Eff-Sc costs slightly more time for $\delta$ around $0.1$ due to the overhead cost of computing $O_\delta$ and $C_\delta$, the computational time stays nearly constant uniformly in $\delta$, indicating that Eff-Sc is a substantially more efficient algorithm than CC-Sc.

\begin{figure}
    \centering
    \begin{subfigure}[b]{0.47\textwidth}
    \centering
    \includegraphics[width=1.1\textwidth]{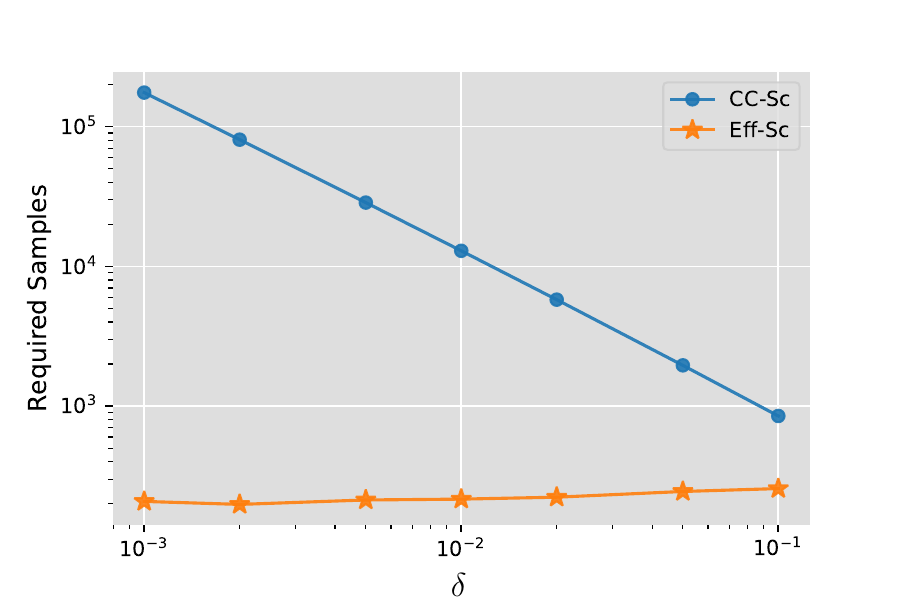}
    \caption{Required Number of Samples in \eqref{eq-portfolio-opt}}
    \label{subfig-portfolio-opt-num-samples}
    \end{subfigure}
    \hfill
    \begin{subfigure}[b]{0.47\textwidth}
    \centering
    \includegraphics[width=1.1\textwidth]{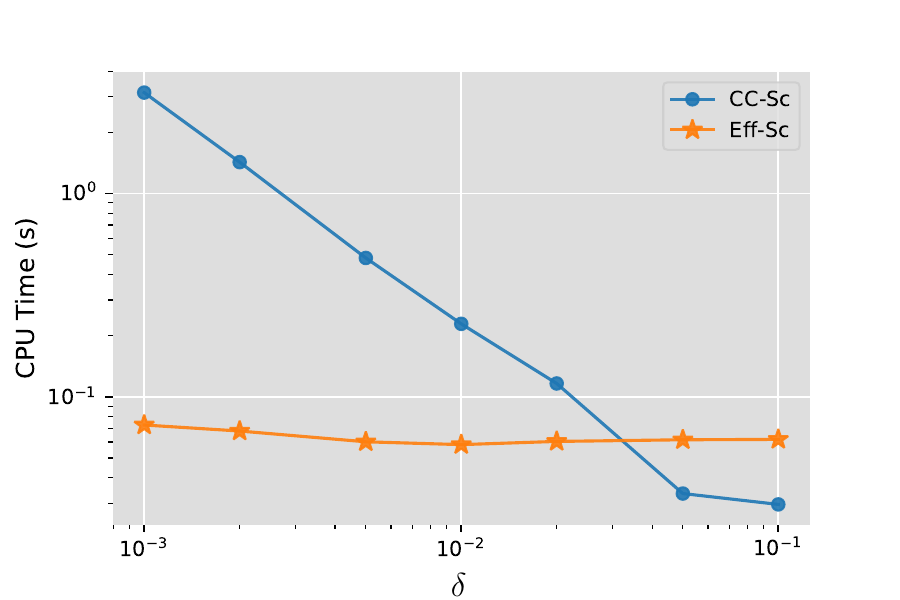}
    \caption{CPU Time}
    \label{subfig-portfolio-opt-cpu-time}
    \end{subfigure}
    \caption{Comparison of computational efficiency for the portfolio optimization problem, in terms of the required number of samples shown in Figure \ref{subfig-portfolio-opt-num-samples} and the used CPU time shown in Figure \ref{subfig-portfolio-opt-cpu-time}. We test  $\delta \in \{0.001,0.002,0.005,0.01,0.02,0.05,0.1\}.$
    }
    \label{fig-efficiency-portfolio-opt}
\end{figure}

Finally, we compare Eff-Sc and CC-Sc for the optimal values of the sampled problems and the violation probabilities of the optimal solutions. Because both methods require generating random samples, the generated solutions are also random. Thus, the optimal values and the violation probabilities are also random. To compare the distributions of the random quantities, we conduct $10^3$ independent experiments. In each experiment, we execute both algorithms and get two solutions, then we evaluate the solutions' violation probabilities using $10^6$ samples of $L$. We employ boxplots (See  \cite{mcgill1978variations}) to depict the samples' distribution through their quantiles. A boxplot is constructed of two parts, a box and a set of whiskers. The box is drawn from the $25\%$ quantile to the $75\%$ quantile, with a horizontal line drawn in the middle to denote the median. Two whiskers indicate $5\%$ and $95\%$ quantiles, respectively, and the scatters represent all the rest sample points beyond the whiskers. 

\begin{figure}
    \centering
    \begin{subfigure}[b]{0.47\textwidth}
    \centering
    \includegraphics[width=1.1\textwidth]{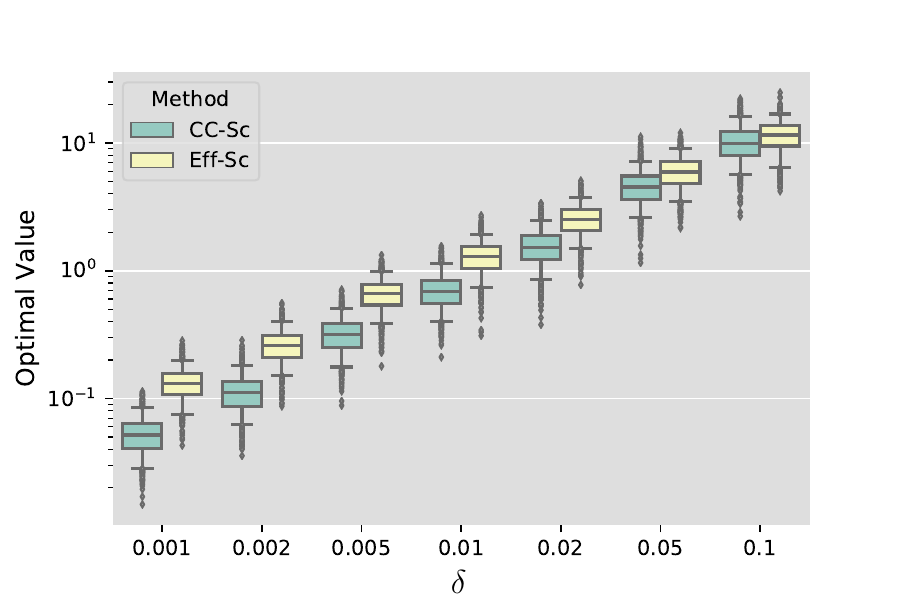}
    \caption{Optimal Value of the Sampled Problem}
    \label{subfig-portfolio-opt-obj-val}
    \end{subfigure}
    \hfill
    \begin{subfigure}[b]{0.47\textwidth}
    \centering
    \includegraphics[width=1.1\textwidth]{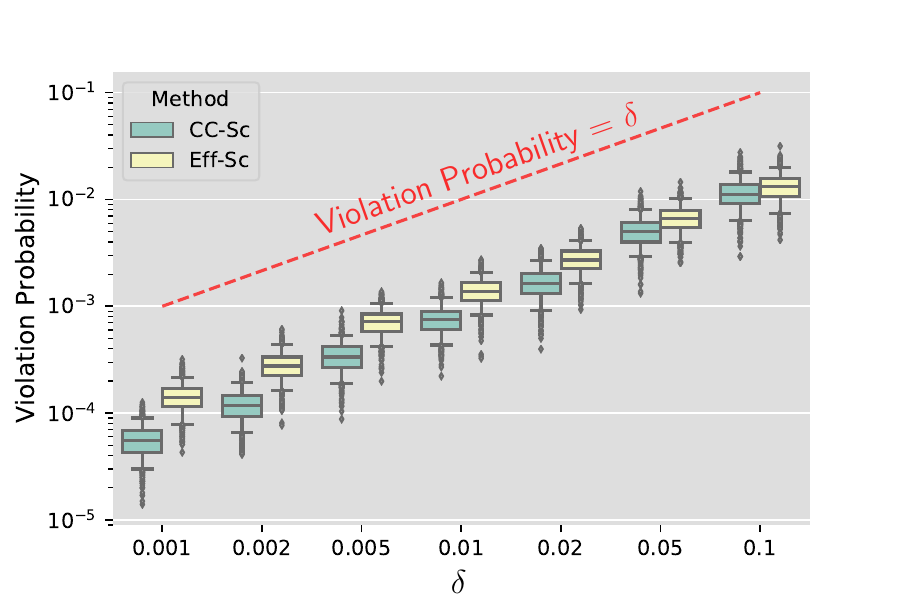}
    \caption{Violation Probability of the Solution}
    \label{subfig-portfolio-opt-prob}
    \end{subfigure}
    \caption{Comparison of the quality of optimal solutions for the portfolio optimization problem, in terms of the optimal value shown in Figure \ref{subfig-portfolio-opt-obj-val} and the solutions' violation probabilities shown in Figure \ref{subfig-portfolio-opt-prob}. Here $\delta \in \{0.001,0.002,0.005,0.01,0.02,0.05,0.1\}$, and the box plots are generated using 1000 experiments.
    }
    \label{fig-solution-quality-portfolio-opt}
\end{figure}

In Figure \ref{fig-solution-quality-portfolio-opt}, we present (a) the optimal values; and (b) the violation probabilities. One can quickly remark from Figure \ref{subfig-portfolio-opt-obj-val} that the optimal value of Eff-Sc is stochastically larger than the optimal value of CC-Sc, while Figure \ref{subfig-portfolio-opt-prob} indicates that the optimal solutions produced by both methods are feasible for all the $10^{3}$ experiments. Overall, with both methods successfully and conservatively approximating the probabilistic constraint, Eff-Sc is more computationally efficient and less conservative, producing solutions with better objective values than its counterpart.

\subsection{Minimal Salvage Fund}
\label{sec-numerics-minimal-salvage-fund}

In this section we conduct a numerical experiment for the minimal salvage
fund problem \eqref{chance-constraint-optimization-general}. In the experiment we pick $d\in\{10,15,20\}$ to test the performance of the problem in different dimensions. 

For each fixed $d$, the parameters of the problem \eqref{chance-constraint-optimization-general} are chosen as follows:
\begin{itemize}
    \item $Q = (Q_{i,j}:i,j
    \in\{1,\cdots,d\})$ where $Q_{i,j} = 1/d$ if $i\neq j$ and otherwise $Q_{i,j}
    = 0$.
    \item $m = (m_{i}:i\in\{1,\cdots,d\})$ where $m_i = 10$ for each $i$.
    \item $L_{i}$ are i.i.d. with Pareto cumulative distribution function $\mathrm{P}(L_i > l) = (1/l)$, for $l \geq 1$.
\end{itemize}

Recall the explicit expressions for sets $O_\delta$ and $C_\delta$ from Proposition \ref{prop-minimal-salvage-fund}. To solve the conditional sampled problem \eqref{conditional-sampled-problem}, it remains to sample $L_{\delta}^{(i)}$ and compute $N'$, the required number of samples. When $\delta$ is small, When $\delta\leq 10^{-3}$, solving the optimization problem \eqref{conditional-sampled-problem} costs much more time than simulating $L^{(i)}_{\delta}$, despite that a simple acceptance rejection scheme is applied to sample $L^{(i)}_{\delta}$ in our experiments. We fix the confidence level parameter $\beta = 10^{-5}$ and set $\delta^{\prime} = \delta/P(L\in C_{\delta}) \geq d^{-1}$, then we can compute $N'$ by the first part of Lemma \ref{lemma-main}. 

Similar to Figure \ref{fig-efficiency-portfolio-opt} of the portfolio optimization problem, we compare the efficiency between Eff-Sc and CC-Sc for different $d$ and $\delta$ in Figure \ref{fig-efficiency-salvage-fund}, in terms of (a) the required number of samples; and (b) the CPU time for solving the sampled approximation problem. We observe that the Eff-Sc has uniformly smaller sample complexity and computational complexity than CC-Sc, where the superiority becomes significant for small $\delta$. In particular, the required number of samples and the used CPU time are bounded for Eff-Sc, while they quickly deteriorate for CC-Sc when $\delta$ becomes smaller. It is also worth noting that Eff-Sc is consistently more efficient than CC-Sc for all the tested dimensions.

\begin{figure}
    \centering
    \begin{subfigure}[b]{0.47\textwidth}
    \centering
    \includegraphics[width=1.1\textwidth]{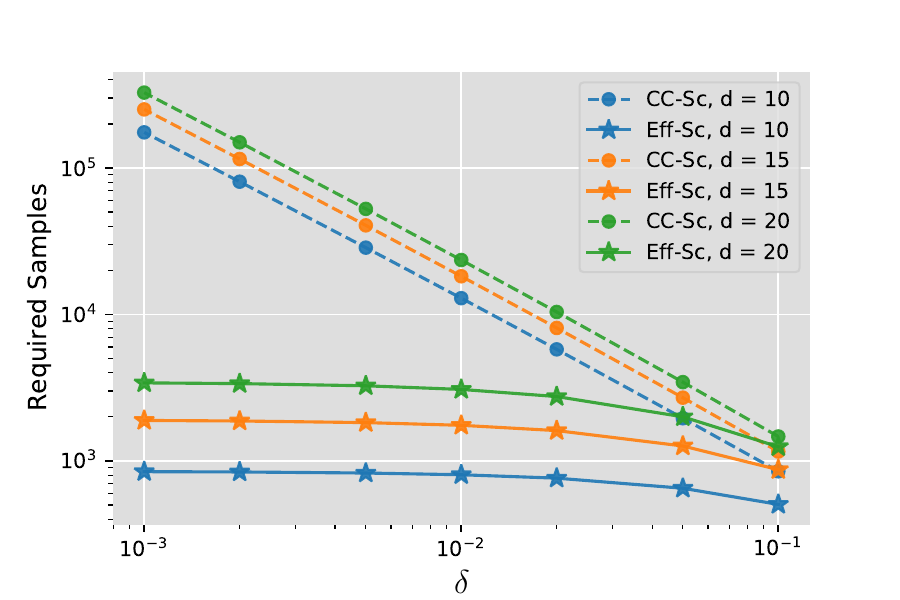}
    \caption{Required Number of Samples in \eqref{chance-constraint-optimization-general}}
    \label{subfig-salvage-fund-num-samples}
    \end{subfigure}
    \hfill
    \begin{subfigure}[b]{0.47\textwidth}
    \centering
    \includegraphics[width=1.1\textwidth]{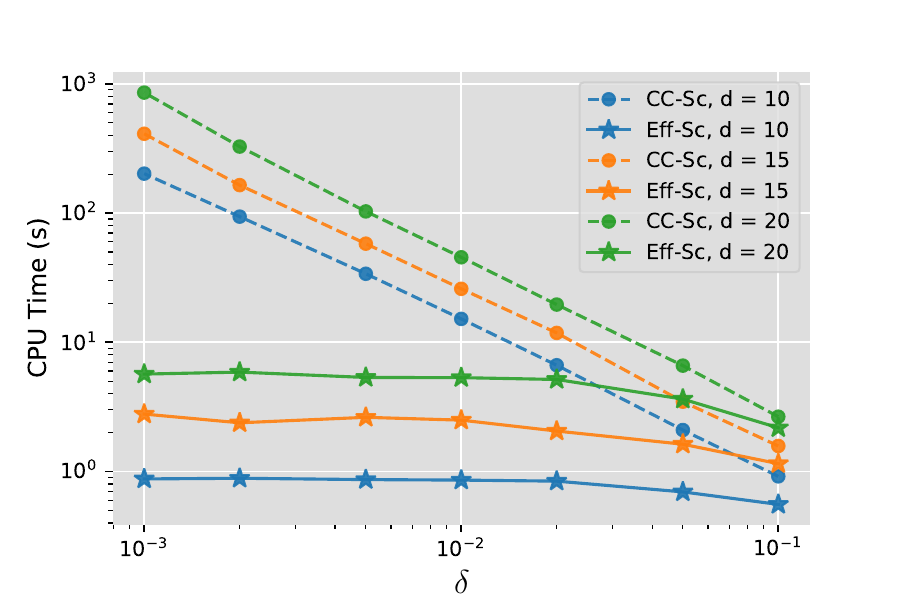}
    \caption{CPU Time}
    \label{subfig-salvage-fund-cpu-time}
    \end{subfigure}
    \caption{Comparison of computational efficiency for the minimal salvage fund problem, in terms of the required number of samples shown in Figure \ref{subfig-salvage-fund-num-samples} and the used CPU time shown in Figure \ref{subfig-salvage-fund-cpu-time}. We test $d\in\{10,15,20\}$ and $\delta \in \{0.001,0.002,0.005,0.01,0.02,0.05,0.1\}.$
    }
    \label{fig-efficiency-salvage-fund}
\end{figure}
\begin{figure}
    \centering
    \begin{subfigure}[b]{0.47\textwidth}
    \centering
    \includegraphics[width=1.1\textwidth]{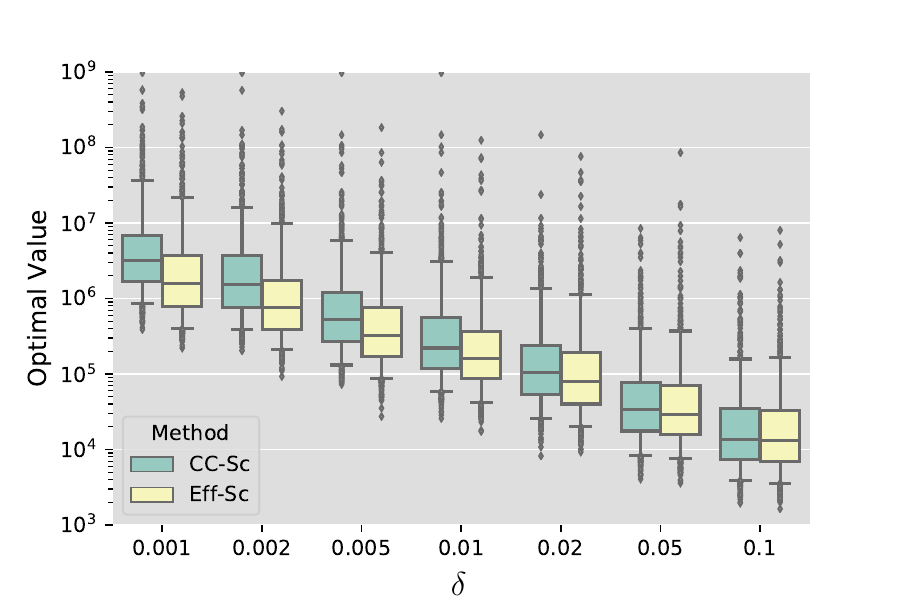}
    \caption{Optimal Value of the Sampled Problem}
    \label{subfig-salvage-fund-obj-val}
    \end{subfigure}
    \hfill
    \begin{subfigure}[b]{0.47\textwidth}
    \centering
    \includegraphics[width=1.1\textwidth]{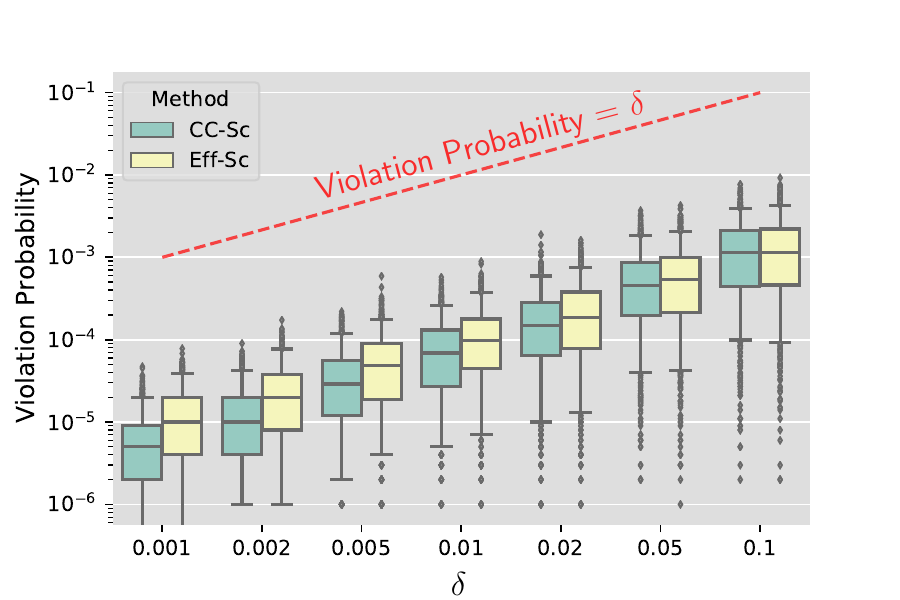}
    \caption{Violation Probability of the Solution}
    \label{subfig-salvage-fund-prob}
    \end{subfigure}
    \caption{Comparison of the quality of optimal solutions for the minimal salvage fund problem, in terms of the optimal value shown in Figure \ref{subfig-salvage-fund-obj-val} and the solutions' violation probabilities shown in Figure \ref{subfig-salvage-fund-prob}. Here $d=15$, $\delta \in \{0.001,0.002,0.005,0.01,0.02,0.05,0.1\}$, and the box plots are generated using 1000 experiments.
    }
    \label{fig-solution-quality-salvage-fund}
\end{figure}
Finally, we compare optimal values of the sampled problems and violation probabilities of the optimal solutions in Figure \ref{fig-solution-quality-salvage-fund}. We present in \eqref{subfig-salvage-fund-obj-val} the optimal values; and \eqref{subfig-salvage-fund-prob} the violation probabilities, with fixed dimension $d = 15$ (We provide additional results for $d = 5$ and $d = 10$ in Appendix \ref{sec-additional-result-minimal-fund}). One can quickly remark from Figure \ref{subfig-salvage-fund-obj-val} that the optimal value of Eff-Sc is stochastically smaller than the optimal value of CC-Sc, while Figure \ref{subfig-salvage-fund-prob} indicates that the optimal solutions produced by both methods are feasible for all the $10^{3}$ experiments. Therefore, we are able to draw the same conclusion as we have from the portfolio optimization experiment: Eff-Sc efficiently produces less conservative solutions. 
\end{revisionenv}

%
%
%


\section*{Acknowledgement}
The authors are grateful to Alexander Shapiro for helpful comments. The research of Bert Zwart is supported by NWO grant 639.033.413. The research of Jose Blanchet is supported by the Air Force Office of Scientific Research under award number FA9550-20-1-0397, NSF grants 1915967, 2118199, 1820942, 1838576, DARPA award N660011824028, and China Merchants Bank. 




\bibliographystyle{informs2014}
\bibliography{ref}

\begin{thebibliography}{42}
\providecommand{\natexlab}[1]{#1}
\providecommand{\url}[1]{\texttt{#1}}
\providecommand{\urlprefix}{URL }

\bibitem[{Ahmed \protect\BIBand{} Shapiro(2008)}]{ahmed2008solving}
Ahmed S, Shapiro A (2008) Solving chance-constrained stochastic programs via
  sampling and integer programming. \emph{State-of-the-Art Decision-Making
  Tools in the Information-Intensive Age}, 261--269 (INFORMS).

\bibitem[{Alsenwi et~al.(2019)Alsenwi, Pandey, Tun, Kim, \protect\BIBand{}
  Hong}]{alsenwi2019chance}
Alsenwi M, Pandey SR, Tun YK, Kim KT, Hong CS (2019) A chance constrained based
  formulation for dynamic multiplexing of embb-urllc traffics in 5g new radio.
  \emph{2019 International Conference on Information Networking (ICOIN)},
  108--113 (IEEE).

\bibitem[{Andrieu et~al.(2010)Andrieu, Henrion, \protect\BIBand{}
  R{\"o}misch}]{andrieu2010model}
Andrieu L, Henrion R, R{\"o}misch W (2010) A model for dynamic chance
  constraints in hydro power reservoir management. \emph{European Journal of
  Operational Research} 207(2):579--589.

\bibitem[{Barrera et~al.(2016)Barrera, Homem-de Mello, Moreno, Pagnoncelli,
  \protect\BIBand{} Canessa}]{barrera2016chance}
Barrera J, Homem-de Mello T, Moreno E, Pagnoncelli BK, Canessa G (2016)
  Chance-constrained problems and rare events: an importance sampling approach.
  \emph{Mathematical Programming} 157(1):153--189.

\bibitem[{Ben-Tal \protect\BIBand{} Nemirovski(2000)}]{ben2000robust}
Ben-Tal A, Nemirovski A (2000) Robust solutions of linear programming problems
  contaminated with uncertain data. \emph{Mathematical programming}
  88(3):411--424.

\bibitem[{Ben-Tal \protect\BIBand{} Nemirovski(2002)}]{ben2002robust}
Ben-Tal A, Nemirovski A (2002) Robust optimization--methodology and
  applications. \emph{Mathematical programming} 92(3):453--480.

\bibitem[{Bertsimas \protect\BIBand{} Sim(2004)}]{bertsimas2004price}
Bertsimas D, Sim M (2004) The price of robustness. \emph{Operations research}
  52(1):35--53.

\bibitem[{Blanchet \protect\BIBand{} Liu(2010)}]{blanchet2010efficient}
Blanchet J, Liu J (2010) Efficient importance sampling in ruin problems for
  multidimensional regularly varying random walks. \emph{Journal of Applied
  Probability} 47(2):301--322.

\bibitem[{Bonami \protect\BIBand{} Lejeune(2009)}]{bonami2009exact}
Bonami P, Lejeune MA (2009) An exact solution approach for portfolio
  optimization problems under stochastic and integer constraints.
  \emph{Operations Research} 57(3):650--670.

\bibitem[{Boyd \protect\BIBand{} Vandenberghe(2004)}]{boyd2004convex}
Boyd S, Vandenberghe L (2004) \emph{Convex optimization} (Cambridge university
  press).

\bibitem[{Calafiore \protect\BIBand{} Campi(2005)}]{calafiore2005uncertain}
Calafiore G, Campi MC (2005) Uncertain convex programs: randomized solutions
  and confidence levels. \emph{Mathematical Programming} 102(1):25--46.

\bibitem[{Calafiore \protect\BIBand{} Campi(2006)}]{calafiore2006scenario}
Calafiore G, Campi MC (2006) The scenario approach to robust control design.
  \emph{IEEE Transactions on Automatic Control} 51(5):742--753.

\bibitem[{Charnes et~al.(1958)Charnes, Cooper, \protect\BIBand{}
  Symonds}]{charnes1958cost}
Charnes A, Cooper WW, Symonds GH (1958) Cost horizons and certainty
  equivalents: an approach to stochastic programming of heating oil.
  \emph{Management Science} 4(3):235--263.

\bibitem[{Chen et~al.(2019)Chen, Blanchet, Rhee, \protect\BIBand{}
  Zwart}]{chen2019efficient}
Chen B, Blanchet J, Rhee CH, Zwart B (2019) Efficient rare-event simulation for
  multiple jump events in regularly varying random walks and compound poisson
  processes. \emph{Mathematics of Operations Research} 44(3):919--942.

\bibitem[{Chen et~al.(2010)Chen, Sim, Sun, \protect\BIBand{}
  Teo}]{chen2010cvar}
Chen W, Sim M, Sun J, Teo CP (2010) From cvar to uncertainty set: Implications
  in joint chance-constrained optimization. \emph{Operations research}
  58(2):470--485.

\bibitem[{Diamond \protect\BIBand{} Boyd(2016)}]{cvxpy}
Diamond S, Boyd S (2016) {CVXPY}: A {P}ython-embedded modeling language for
  convex optimization. \emph{Journal of Machine Learning Research} 17(83):1--5.

\bibitem[{Eisenberg \protect\BIBand{} Noe(2001)}]{eisenberg2001systemic}
Eisenberg L, Noe TH (2001) Systemic risk in financial systems. \emph{Management
  Science} 47(2):236--249.

\bibitem[{Embrechts et~al.(2013)Embrechts, Kl{\"u}ppelberg, \protect\BIBand{}
  Mikosch}]{embrechts2013modelling}
Embrechts P, Kl{\"u}ppelberg C, Mikosch T (2013) \emph{Modelling extremal
  events: for insurance and finance}, volume~33 (Springer Science \& Business
  Media).

\bibitem[{Frank(2008)}]{frank2008municipal}
Frank B (2008) Municipal bond fairness act. \emph{110th Congress, 2d Session,
  House of Representatives, Report}, 110--835.

\bibitem[{Gudmundsson \protect\BIBand{} Hult(2014)}]{gudmundsson2014markov}
Gudmundsson T, Hult H (2014) Markov chain monte carlo for computing rare-event
  probabilities for a heavy-tailed random walk. \emph{Journal of Applied
  Probability} 51(2):359--376.

\bibitem[{Hillier(1967)}]{hillier1967chance}
Hillier FS (1967) Chance-constrained programming with 0-1 or bounded continuous
  decision variables. \emph{Management Science} 14(1):34--57.

\bibitem[{Hong et~al.(2020)Hong, Huang, \protect\BIBand{}
  Lam}]{hong2020learning}
Hong LJ, Huang Z, Lam H (2020) Learning-based robust optimization: Procedures
  and statistical guarantees. \emph{Management Science} .

\bibitem[{Hong et~al.(2011)Hong, Yang, \protect\BIBand{}
  Zhang}]{hong2011sequential}
Hong LJ, Yang Y, Zhang L (2011) Sequential convex approximations to joint
  chance constrained programs: A monte carlo approach. \emph{Operations
  Research} 59(3):617--630.

\bibitem[{Kley et~al.(2016)Kley, Kl{\"u}ppelberg, \protect\BIBand{}
  Reinert}]{kley2016risk}
Kley O, Kl{\"u}ppelberg C, Reinert G (2016) Risk in a large claims insurance
  market with bipartite graph structure. \emph{Operations Research}
  64(5):1159--1176.

\bibitem[{K{\"u}{\c{c}}{\"u}kyavuz(2012)}]{kuccukyavuz2012mixing}
K{\"u}{\c{c}}{\"u}kyavuz S (2012) On mixing sets arising in chance-constrained
  programming. \emph{Mathematical programming} 132(1-2):31--56.

\bibitem[{Lagoa et~al.(2005)Lagoa, Li, \protect\BIBand{}
  Sznaier}]{lagoa2005probabilistically}
Lagoa CM, Li X, Sznaier M (2005) Probabilistically constrained linear programs
  and risk-adjusted controller design. \emph{SIAM Journal on Optimization}
  15(3):938--951.

\bibitem[{Lejeune \protect\BIBand{} Margot(2016)}]{lejeune2016solving}
Lejeune MA, Margot F (2016) Solving chance-constrained optimization problems
  with stochastic quadratic inequalities. \emph{Operations Research}
  64(4):939--957.

\bibitem[{Luedtke(2014)}]{luedtke2014branch}
Luedtke J (2014) A branch-and-cut decomposition algorithm for solving
  chance-constrained mathematical programs with finite support.
  \emph{Mathematical Programming} 146(1-2):219--244.

\bibitem[{Luedtke \protect\BIBand{} Ahmed(2008)}]{luedtke2008sample}
Luedtke J, Ahmed S (2008) A sample approximation approach for optimization with
  probabilistic constraints. \emph{SIAM Journal on Optimization}
  19(2):674--699.

\bibitem[{Luedtke et~al.(2010)Luedtke, Ahmed, \protect\BIBand{}
  Nemhauser}]{luedtke2010integer}
Luedtke J, Ahmed S, Nemhauser GL (2010) An integer programming approach for
  linear programs with probabilistic constraints. \emph{Mathematical
  Programming} 122(2):247--272.

\bibitem[{McGill et~al.(1978)McGill, Tukey, \protect\BIBand{}
  Larsen}]{mcgill1978variations}
McGill R, Tukey JW, Larsen WA (1978) Variations of box plots. \emph{The
  American Statistician} 32(1):12--16.

\bibitem[{{MOSEK ApS}(2020)}]{mosek}
{MOSEK ApS} (2020) \emph{MOSEK Fusion API for Python}.
  \urlprefix\url{https://docs.mosek.com/9.2/pythonfusion.pdf}.

\bibitem[{Nemirovski \protect\BIBand{}
  Shapiro(2006{\natexlab{a}})}]{nemirovski2006convex}
Nemirovski A, Shapiro A (2006{\natexlab{a}}) Convex approximations of chance
  constrained programs. \emph{SIAM Journal on Optimization} 17(4):969--996.

\bibitem[{Nemirovski \protect\BIBand{}
  Shapiro(2006{\natexlab{b}})}]{nemirovski2006scenario}
Nemirovski A, Shapiro A (2006{\natexlab{b}}) Scenario approximations of chance
  constraints. \emph{Probabilistic and Randomized Methods for Design under
  Uncertainty}, 3--47 (Springer).

\bibitem[{Pe{\~n}a-Ordieres et~al.(2020)Pe{\~n}a-Ordieres, Luedtke,
  \protect\BIBand{} Wächter}]{pena2020solving}
Pe{\~n}a-Ordieres A, Luedtke JR, Wächter A (2020) Solving chance-constrained
  problems via a smooth sample-based nonlinear approximation. \emph{SIAM
  Journal on Optimization} 30(3):2221--2250.

\bibitem[{Prekopa(1970)}]{prekopa1970probabilistic}
Prekopa A (1970) On probabilistic constrained programming. \emph{Proceedings of
  the Princeton Symposium on Mathematical Programming}, volume 113, 138
  (Princeton, NJ).

\bibitem[{Pr{\'e}kopa(2003)}]{prekopa2003probabilistic}
Pr{\'e}kopa A (2003) Probabilistic programming. \emph{Handbooks in Operations
  Research and Management Science} 10:267--351.

\bibitem[{Resnick(2013)}]{resnick2013extreme}
Resnick SI (2013) \emph{Extreme values, regular variation and point processes}
  (Springer).

\bibitem[{Sepp{\"a}l{\"a}(1971)}]{seppala1971constructing}
Sepp{\"a}l{\"a} Y (1971) Constructing sets of uniformly tighter linear
  approximations for a chance constraint. \emph{Management Science}
  17(11):736--749.

\bibitem[{Tong et~al.(2020)Tong, Subramanyam, \protect\BIBand{}
  Rao}]{tong2020optimization}
Tong S, Subramanyam A, Rao V (2020) Optimization under rare chance constraints.
  \emph{arXiv preprint arXiv:2011.06052} .

\bibitem[{Wierman \protect\BIBand{} Zwart(2012)}]{wierman2012tail}
Wierman A, Zwart B (2012) Is tail-optimal scheduling possible? \emph{Operations
  Research} 60(5):1249--1257.

\bibitem[{Zhang et~al.(2014)Zhang, K{\"u}{\c{c}}{\"u}kyavuz, \protect\BIBand{}
  Goel}]{zhang2014branch}
Zhang M, K{\"u}{\c{c}}{\"u}kyavuz S, Goel S (2014) A branch-and-cut method for
  dynamic decision making under joint chance constraints. \emph{Management
  Science} 60(5):1317--1333.

\end{thebibliography}

\clearpage
\appendix
\section{Proofs of Technical Results}\label{sec-proofs}
\subsection{Proofs for Section \ref{sec-construction}}

\begin{proofenv}{Proof of Lemma \ref{lemma-outer-approx}}
    We will derive an expression of $\alpha_{\delta}$ to ensure that $F_{\delta
	}\subseteq\bigcup_{\alpha\geq\alpha_{\delta}}\alpha\cdot\Pi$ for $\delta$
	small enough. Because of Assumption \ref{assumption-level-set}, for any
	$\alpha_{0} >0$ there exist some $\delta$ small enough such that $F_{\delta
	}\subseteq\bigcup_{\alpha\geq\alpha_{0}} \alpha\cdot\Pi$. Therefore, it
	suffices to prove that $F_{\delta}$ and $\bigcup_{\alpha< \alpha_{\delta}%
	}\alpha\cdot\Pi$ are disjoint. In other words,
	\begin{align}
	\label{eq-to-prove}\mathrm{P}\left( \phi(\alpha\cdot x,L)>0\right)  >
	\delta,\quad\forall\alpha<\alpha_{\delta}, x\in\Pi, \delta<\delta_{0}.
	\end{align}
	
	Let $\varepsilon$ be a positive number such that $\mu(C_{\varepsilon, -}) >
	0$. Pick the set $K$ in \eqref{eq-asymp-lower} as a compact set such that
	$0<\mu(K\cap C_{\varepsilon, -})<\infty$. It follows from the inequality
	\eqref{eq-asymp-lower} that there exist a constant $\alpha_{1}$ such that
	\begin{equation}
	\label{eq-bounded-diff-min}\Psi_{-}(l) - \varepsilon\leq\inf_{x\in\Pi}\left[
	\frac{1}{r(\alpha)}\phi(\alpha\cdot x, h(\alpha)\cdot l)\right]  \qquad\forall
	l\in K,\; \alpha> \alpha_{1}%
	\end{equation}
	Therefore, for any $\alpha\geq\alpha_{1}$ we have,
	\begin{align}
	\label{eq-lower-prob-bound}%
	\begin{split}
	\mathrm{P}\left( \min_{x\in\Pi}\phi(\alpha\cdot x,L)>0\right)   & =
	\mathrm{P}\left( \min_{x\in\Pi}\frac{1}{r(\alpha)}\phi(\alpha\cdot
	x,L)>0\right) \\
	(\text{Due to }\eqref{eq-bounded-diff-min}) & \geq\mathrm{P}\left(
	G(L/h(\alpha))\geq\varepsilon; L/h(\alpha)\in K\right) \\
	& = \mathrm{P}\left( L\in h(\alpha)\cdot(K\cap C_{\varepsilon, -}) \right) .
	\end{split}
	\end{align}
	Recall that $L$ is regularly varying from Assumption
	\ref{assumption-regular-varying},
	\begin{align*}
	\lim_{\alpha\rightarrow\infty} \frac{\mathrm{P}\left( L\in h(\alpha
		)\cdot(K\cap C_{\varepsilon, -})\right) }{\mathrm{P}(\|L\|_{2}>h(\alpha))} =
	\mu(K\cap C_{\varepsilon, -}).
	\end{align*}
	Therefore, there exist a number $\alpha_{2}$ such that
	\begin{align}
	\label{eq-asymp-prob}\mathrm{P}\left( L\in h(\alpha)\cdot(K\cap C_{\varepsilon
		, -})\right)  \geq\frac{1}{2} \mathrm{P}(\|L\|_{2}>h(\alpha))\mu(K\cap
	C_{\varepsilon, -}), \qquad\forall\alpha\geq\alpha_{2}.
	\end{align}
	Note that the right hand side of \eqref{eq-asymp-prob} is nondecreasing in
	$\alpha$. Thus, if $\delta_{1}:=\frac{1}{2} \mathrm{P}\big(\|L\|_{2}%
	>h(\alpha_{2})\big)\mu(K\cap C_{\varepsilon, -}),$  for any $\delta\leq
	\delta_{1}$ there exist $\alpha_{\delta}$ satisfying
	\begin{align}
	\label{eq-b-delta}\frac{1}{2} \mathrm{P}(\|L\|_{2}>h(\alpha_{\delta}%
	))\mu(K\cap C_{\varepsilon, -}) = \delta.\qquad\forall\alpha,\delta\quad
	s.t.\quad\alpha_{2}\leq\alpha<\alpha_{\delta}, 0<\delta\leq\delta_{1}.
	\end{align}
	Substituting \eqref{eq-b-delta} into \eqref{eq-lower-prob-bound}, we have
	\begin{align*}
	& \mathrm{P}\left( \phi(x,L)>0\right) \geq\mathrm{P}\left( \min_{x\in\Pi}%
	\phi(\alpha\cdot x,L)>0\right)  > \delta.\\
	& \forall\alpha,x,\delta\quad s.t.\quad\max(\alpha_{1},\alpha_{2})\leq
	\alpha<\alpha_{\delta}, x\in\Pi, 0<\delta\leq\delta_{1}.
	\end{align*}
	Moreover, Assumption \ref{assumption-level-set} guarantees the existence of
	$\delta_{2}$ such that
	\begin{align*}
	\mathrm{P}\left( \phi(\alpha\cdot x,L)>0\right)  > \delta,\quad\forall
	\alpha<\max(\alpha_{1},\alpha_{2}), x\in\Pi, \delta<\delta_{2}.
	\end{align*}
	Consequently \eqref{eq-to-prove} is proved with $\delta_{0} = \min(\delta
	_{1},\delta_{2})$.
\end{proofenv}

\begin{proofenv}{Proof of Theorem \ref{thm-set-construction}}
We construct the uniform conditional event $C_{\delta}$ that
	contains all the $V_{x}$ for $x\in O_{\delta}$. Due to the definition
	\eqref{eq-asymp-upper} and $\lim_{\delta\rightarrow0} \alpha_{\delta}= \infty
	$, there exist $\delta_{0}$ such that for all $\delta<\delta_{0}$,
	\begin{equation}
	\label{eq-bounded-diff-max}\Psi_{+}(l) + \varepsilon\geq\sup_{x\in\Pi}\left[
	\frac{1}{r(\alpha)}\phi(\alpha\cdot x, h(\alpha)\cdot l)\right]  \qquad\forall
	l\in K,\; \alpha> \alpha_{\delta}.%
	\end{equation}
	Notice that for any $x\in O_{\delta}$, there exist an $\alpha_{x}\geq
	\alpha_{\delta}$ such that $x\in\alpha_{x}\cdot\Pi$. Consequently, it follows
	from \eqref{eq-bounded-diff-max} that
	\begin{align*}
	\phi(x,l)>0 \Longrightarrow\Psi_{+}\Big(\frac{l}{h(\alpha_{x})}\Big) \geq
	-\varepsilon, \qquad\forall x\in O_{\delta}, l\in h(\alpha_{x})\cdot K.
	\end{align*}
	Applying Assumption \ref{assumption-approx-func} yields that
	\begin{align*}
	\Psi_{+}\Big(\frac{l}{h(\alpha_{\delta})}\Big) \geq\Psi_{+}\Big(\frac
	{l}{h(\alpha_{x})}\Big) \geq-\varepsilon, \qquad\forall x\in O_{\delta}, l\in
	h(\alpha_{x})\cdot(K\cap S).
	\end{align*}
	Recall that $K$ is a ball in $\mathbb{R}^{d_{l}}$ (thus $K\subseteq
	(h(\alpha_{x})/h(\alpha_{\alpha}))\cdot K$) and that $S\subseteq(h(\alpha
	_{x})/h(\alpha_{\alpha}))\cdot S$ from Assumption \ref{assumption-approx-func}%
	, it turns out that $h(\alpha_{\delta})\cdot(K\cap S)\subseteq h(\alpha
	_{x})\cdot(K\cap S)$.  Consequently, whenever $l\in V_{x}$ for some $x\in
	O_{\delta}$, we either have $l \in h(\alpha_{x})\cdot(K\cap S)$ implying
	$\Psi_{+}\Big(\frac{l}{h(\alpha_{\delta})}\Big) \geq-\varepsilon$, or we have
	$l\in\big(h(\alpha_{x})\cdot(K\cap S)\big)^{c}\subseteq(h(\alpha_{\delta
	})\cdot(K\cap S)\big)^{c}$. Summarizing these two scenarios,
	\begin{align*}
	\bigcup_{x\in O_{\delta}} V_{x}  & \subseteq\{l\in\mathbb{R}^{d_{l}}\mid
	\Psi_{+}\Big(\frac{l}{h(\alpha_{\delta})}\Big) \geq-\varepsilon\}\bigcup
	\big(h(\alpha_{\delta})\cdot(K\cap S)\big)^{c}\\
	&  = h(\alpha_{\delta})\cdot\big(C_{\varepsilon,+} \cup K^{c}\cup S^{c}\big).
	\end{align*}
	Thus, we define the conditional set $C_{\delta}$ as
	\begin{align*}
	C_{\delta} & := h(\alpha_{\delta})\cdot\big(C_{\varepsilon,+} \cup K^{c}\cup
	S^{c}\big).
	\end{align*}

	It remains to analyze the probability of the uniform conditional event
	$C_{\delta}$. As $L$ is multivariate regularly varying,
	\begin{align*}
	\lim_{\delta\rightarrow0} \frac{\mathrm{P}\left( L\in C_{\delta}\right)
	}{\mathrm{P}(\|L\|_{2}>h(\alpha_{\delta}))} = \mu(C_{\varepsilon,+} \cup
	K^{c}\cup S^{c}).
	\end{align*}
	Recalling,  $\mathrm{P}(\|L\|_{2}>h(\alpha_{\delta})) = O(\delta)$ and invoking Property
	\ref{condition-outer-approx}, we get
	\begin{align*}
	\limsup_{\delta\rightarrow0}\delta^{-1}\mathrm{P}(L\in C_{\delta}) < \infty.
	\end{align*}
	Hence, the proof is complete.
\end{proofenv}

\begin{proofenv}{Proof of Theorem \ref{thm-nearly-optimal}}
Using Lemma \ref{lemma-main}, we immediately have 
	$ \mathrm{P}(\Val	\eqref{chance-constraint-opt} \leq \Val\eqref{conditional-sampled-problem})
		\geq 1-\beta,$
	it remains to show that there exist $\Lambda>0$ such that 
	$ \mathrm{P}(\Val \eqref{conditional-sampled-problem}
	\leq \Lambda \times \Val \eqref{chance-constraint-opt})
	\geq 1-\beta$.
	
	For simplicity, in the proof we will use $L_\delta$ as a shorthand for $(L|L\in C_\delta)$, the random variable with conditional distribution of $L$ given $L\in C_\delta$. By a scaling of $x$ by a factor $\alpha_\delta$ in \eqref{conditional-sampled-problem}, we have an equivalent optimization problem
	\begin{align}
	\label{scaled-conditional-problem}
	\begin{array}
	[c]{ll}%
	\mbox{minimize} & c^{\top}x\\
	\mbox{subject to} & \frac{1}{r(\alpha_\delta)}\phi(\alpha_\delta\cdot x, L_{\delta}^{(i)}) \leq0, \quad i = 1,\ldots
	,N^{\prime},\\
	& x\in \bigcup_{\alpha\geq1}\alpha
	\cdot\Pi.
	\end{array}
	\end{align}
	where $L_{\delta}^{(i)}$ are i.i.d. samples from $L_\delta$. Notice that  $\Val\eqref{conditional-sampled-problem} = \alpha_\delta \times \Val\eqref{scaled-conditional-problem}$.
	
	For any compact set $E\subseteq C_\delta$, since $L$ is multivariate regularly varying,
	\[
	\lim_{\delta\rightarrow 0} 
	\mathrm{P} ((h(\alpha_\delta))^{-1} L_\delta \in E)
	=  \lim_{\delta\rightarrow 0} 
	\frac{\mathrm{P} (L \in(h(\alpha_\delta)\cdot E))}
	{\mathrm{P}(L\in C_\delta)}
	= \frac{\lim_{\delta\rightarrow 0}
		\frac{\mathrm{P} (L \in(h(\alpha_\delta)\cdot E))}{\mathrm{P}(\|L\|_{2}>h(\alpha_{\delta}))}
	}{\lim_{\delta\rightarrow 0}
		\frac{\mathrm{P} (L \in C_\delta)}{\mathrm{P}(\|L\|_{2}>h(\alpha_{\delta}))}}
	= \frac{\mu(E)}{\mu(C_{\varepsilon, +}
		\cup K^{c}\cup S^{c})}. 
	\]
	Thus $(h(\alpha_\delta))^{-1} L_\delta \vconverge L_{\lim}$. As the limiting measure is a probability measure, the family $\{h(\alpha_\delta))^{-1} L_\delta \mid \delta > 0\}$ is tight and consequently $(h(\alpha_\delta))^{-1} L_\delta \dconverge L_{\lim}$ follows directly from the vague convergence, see \cite{resnick2013extreme}. Consequently, since all the samples are i.i.d, we also have 
	\[
	(h(\alpha_\delta))^{-1} \cdot (L_{\delta}^{(1)}, \ldots, L_{\delta}^{(N')})
	\dconverge (L_{\lim}^{(1)}, \ldots, L_{\lim}^{(N')}).
	\]
	
	Now we define a family of deterministic optimization problem, denoted by \eqref{deterministic-problem}, which is parameterized by $(l_1,\cdots, l_{N'})$ as follows, 
	\begin{align}
	\begin{array}
	[c]{ll}%
	\mbox{minimize} & c^{\top}x\\
	\mbox{subject to} & \phi_{\lim}(x, l_i) \leq0, \quad i = 1,\ldots
	,N^{\prime},\\
	& x\in \bigcup_{\alpha\geq1}\alpha
	\cdot\Pi.
	\end{array}
	\label{deterministic-problem}
	\tag{$DP(l_1,\cdots, l_{N'})$}
	\end{align}
	Then, there exist a compact set $E_1 \subseteq \mathbb{R}^{d_l\times N'}$ such that:
	\begin{enumerate}
		\item Problem \eqref{deterministic-problem} satisfies Slater's condition if $(l_1,\cdots, l_{N'})\in E_1$;
		\item $\mathrm{P}((h(\alpha_\delta))^{-1} \cdot (L_{\delta}^{(1)}, \ldots, L_{\delta}^{(N')}) \in E_1) \geq 1-\beta$ for all $\delta>0$;
	\end{enumerate}

  	For every $(l_1,\cdots, l_{N'})\in E_1$ and $\epsilon>0$, due to the Slater's condition, there exist a feasible solution $x\in \bigcup_{\alpha\geq1}\alpha$ such that $\sup_{j = 1,\ldots,N'}  		\phi_{\lim}(x, l_j) < -\epsilon$. Since $\phi_{\lim}(x,l)$ is continuous in $l$, there exist an open neighborhood $U$ around $(l_1,\cdots, l_{N'})$ such that $\sup_{(l_1,\ldots, l_{N'})\in U}\sup_{j = 1,\ldots, N'} \phi_{\lim}(x, l_j) < -\epsilon/2$. Notice that such feasible solution $x$ and neighborhood $U$ exist for every $(l_1,\cdots, l_{N'})\in E_1$. 
  	\begin{revisionenv}
  	There exists a finite open cover $\{U_i\}_{i=1}^m$ of $E_1$ due to its compactness. Let $\{x_i\}_{i=1}^m$ be the corresponding feasible solutions to the open cover $\{U_i\}_{i=1}^m$.
  	  \end{revisionenv}
  	Due to Assumption \ref{assumption-lim-phi}, there exist $\delta_1 > 0$ such that for all $\delta<\delta_1$, we have
  	\begin{align}
  	\sup_{(l_1,\ldots, l_{N'})\in E_1}
  	\sup_{i = 1,\ldots,m}
  	\sup_{j = 1,\ldots,N'}  		
  	\left\vert
  	\frac{1}{r(\alpha_\delta)}\phi(\alpha_\delta\cdot x_i, h(\alpha_\delta)\cdot l_j)
  	- \phi_{\lim}(x_i,l_j)
  	\right\vert < \epsilon/2. 
  	\label{eq-approximation-phi-lim}
  	\end{align} 
  	Therefore by the triangle inequality, it follows that if $\delta<\delta_1$, 
  	\begin{align*}
  	\sup_{(l_1,\ldots, l_{N'})\in U_i}\sup_{j = 1,\ldots, N'} \frac{1}{r(\alpha_\delta)}\phi(\alpha_\delta\cdot x_i, h(\alpha_\delta)\cdot l_j) < 0 .
  	\end{align*}
  	Consequently, $x_i$ is a feasible solution for optimization problem \eqref{scaled-conditional-problem} if $(h(\alpha_\delta))^{-1} \cdot(L_{\delta}^{(1)},\ldots, L_{\delta}^{(N')})\in U_i$, which further implies that $\alpha_\delta^{-1}\times\Val\eqref{conditional-sampled-problem}\leq c^{\top}x_i$. As a result, we have 
  	$$\Val\eqref{conditional-sampled-problem}\leq \alpha_\delta \times \max_{i = 1,\ldots,m}c^{\top}x_i,
  	 \quad\mbox{if}\quad (h(\alpha_\delta))^{-1} \cdot(L_{\delta}^{(i)},\ldots, L_{\delta}^{(N')})\in E_1.$$ Note that $\Val\eqref{chance-constraint-opt}\geq \inf_{x\in O_\delta}c^{\top}x = \alpha_\delta \times 
  	 \inf\{c^{\top}x\mid x\in \bigcup_{\alpha\geq1}\alpha
  	 \cdot\Pi\}$. Therefore, let 
  	 $$\Lambda = \Big(\inf\{c^{\top}x\mid x\in \bigcup_{\alpha\geq1}\alpha\cdot\Pi\}\Big)^{-1}
  	 \times\big(\max_{i = 1,\ldots,m}c^{\top}x_i \big)>0$$ 
  	 It follows that $$ \mathrm{P}\Big(\Val \eqref{conditional-sampled-problem}
  	 \leq \Lambda \times \Val \eqref{chance-constraint-opt}\Big)
  	 \geq \mathrm{P}\Big((h(\alpha_\delta))^{-1} \cdot(L_{\delta}^{(1)},\ldots, L_{\delta}^{(N')})\in E_1\Big)
  	 \geq 1-\beta.$$
  	 The statement is concluded by using the union bound, combining the lower bound together with the upper bound implied by Theorems \ref{lemma-main} and \ref{thm-set-construction}, hence obtaining factor $2\beta$.
\end{proofenv}

\begin{revisionenv}
\begin{proofenv}{Proof of Lemma \ref{lemma-linear-approx}.}
	Without loss of generality, assume that $R$ is an integer such that
	\begin{align*}
	Z_{\phi}= \{(x,l)\in\mathbb{R}^{d_{x}}\times\mathbb{R}^{d_{l}}\mid
	\phi(x,l)\leq0 \} \subseteq[-R,R]^{(d_{x}+d_{l})}.
	\end{align*}
	Let $N_1 = (2R+1)^{(d_{x}+d_{l})}$, and let $(x^{(i)},l^{(i)}), i = 1,\ldots,N_1$ be the integer lattice points in $[-R,R]^{(d_{x}+d_{l})}$. In addition, let
	$a_{i} = \frac{\partial\phi}{\partial L}(x^{(i)},l^{(i)})$, $b_{i} =
	\frac{\partial\phi}{\partial x}(x^{(i)},l^{(i)})$ and $c_{i} = \phi(x^{(i)},l^{(i)})
	- \frac{\partial\phi}{\partial L}(x^{(i)},l^{(i)})^{\top}l^{(i)} - \frac{\partial\phi
	}{\partial x}(x^{(i)},l^{(i)})^{\top}x^{(i)}$ for $i = 1,\ldots,N_1$, then define $\phi_{1,-}(x, l) = \max_{i = 1,\ldots, N_1}a_{i}^{\top}l + b_{i}^{\top} x + c_{i}$. Since the function  $\phi(x, l)$ is convex, we can invoke the supporting hyperplane theorem to deduce that  $a_{i}^{\top}l + b_{i}^{\top} x + c_{i}\leq\phi(x,l)$ for $i = 1,\ldots, N_1$, and
	consequently $\phi_{1, -}(x,l)\leq\phi(x,l)$. In addition, since $\phi(x,l) \geq 0$ at the boundary of the cube $[-R,R]^{(d_{x}+d_{l})}$, there exist a constant $C_1$ such that $-C_1\cdot R \pm C_1\cdot x_i\leq \phi(x,l)$ for $i = 1,\ldots, d_x$ and $-C_1\cdot R \pm C_1\cdot l_i\leq \phi(x,l)$ for $i = 1,\ldots, d_l$, for all $(x,l)\in\mathbb{R}^{d_x} \times\mathbb{R}^{d_l}$. Therefore, with $\phi_{2, -}(x,l)$ being the maximum of the aforementioned $N_2 = 2(d_x+d_l)$ linear functions, we have $\phi_{2, -}(x,l)\leq \phi(x,l)$, and we also have that $\phi_{2, -}(x,l) \leq 0$ implies $(x,l)\in [-R,R]^{(d_{x}+d_{l})}$.
	
	Define $\phi_{-}(x,l) = \max\left\{\phi_{1, -}(x,l),\phi_{2, -}(x,l)\right \}$. We can conclude the property of $\phi_{-}(x,l)$ as follows: (1) $\phi_{-}(x,l)$ is a piecewise linear function of form  $\max_{i = 1,\ldots, N}a_{i}^{\top}l + b_{i}^{\top} x + c_{i}$, where $N = N_1 + N_2$; (2) $\phi_{-}(x,l)\leq \phi(x,l)$; (3) $\phi_{-}(x,l) \leq 0$ implies $(x,l)\in [-R,R]^{(d_{x}+d_{l})}$. To complete the proof, it remains to verify for $\phi_{-}(x,l)$ the second statement of Assumption \ref{assumption-linear-approx}.
	
	As $\phi_{-}(x,l)\leq 0$ implies $(x,l)\in [-R,R]^{(d_{x}+d_{l})}$, it suffices to prove that there exist some universal constant $C\in\mathbb{R}_{+}$ such that $\phi(x,l) - \phi_{-}(x,l)\leq C$ for all $(x,l)\in [-R,R]^{(d_{x}+d_{l})}$. For an arbitrary point $(x,l)\in [-R,R]^{(d_{x}+d_{l})} $, there exist a lattice point $(x^{(i)},l^{(i)})$ such that
	$\|(x,l)-(x^{(i)},l^{(i)})\|_{2}\leq\sqrt{d_{x}+d_{l}}/2.$ Next, since
	$\phi(x,l)$ is twice continuously differentiable, the gradient $\nabla
	\phi(x,l)$ is Lipschitz over $[-R,R]^{(d_{x}+d_{l})}$ with Lipschitz constant denoted by
	$M_{\phi}$. Therefore, for any $(x,l)\in [-R,R]^{(d_{x}+d_{l})}$,
	\begin{align*}
	\phi(x,l) - \phi_{-}(x,l)\leq\phi(x,l) - \phi_{1,-}(x,l)\leq \min_{i = 1,\ldots,N_1}\big(\phi(x,L) - (a_{i}^{\top}L +
	b_{i}^{\top} x + c_{i})\big)\leq\frac{1}{4}M^{2}_{\phi}\sqrt{d_{x}+d_{l}}.
	\end{align*}
	The proof is now complete.
\end{proofenv}

\begin{proofenv}{Proof of Theorem \ref{thm-convex}.}
	Since $\phi_{-}(x,L)\leq\phi(x,L)$, the probability constraint $\mathrm{P}%
	(\phi(x,L)>0)\leq\delta$ implies that $\mathrm{P}(\phi_{-}(x,L)>0)\leq\delta$,
	which further implies $\mathrm{P}(a_{i}^{\top}L + b_{i}^{\top} x + c_{i}>0)\leq\delta$ for $i
	= 1,\ldots,N$. Therefore, we have $-b_{i}^{\top}x - c_{i} \geq\bar{F}^{-1}%
	_{a_{i}^{\top}L}(\delta)$ for $i = 1,\ldots,N$, which implies $F_{\delta}\subseteq
	O_{\delta}$.
	
	Then, consider $x\in O_{\delta}$ and $L\in V_{x} = \{L\in\mathbb{R}^{d_{l}%
	}\mid\phi(x,L)>0\}$. It follows from the second statement of Assumption \ref{assumption-linear-approx} that
	$\phi(x, L)>0$ implies that $\phi_{-}(x, L)+C>0$. Thus, there exist an index
	$i$ such that $a_{i}^{\top}L + b_{i}^{\top} x + c_{i} + C>0$. As $x\in O_{\delta}$
	implies that $b_{i}^{\top}x + c_{i} + \bar{F}^{-1}_{a_{i}^{\top}L}(\delta)\leq0$,
	so
	\begin{align*}
	a_{i}^{\top}L - \bar{F}^{-1}_{a_{i}^{\top}L}(\delta)+C \geq a_{i}^{\top}L + b_{i}^{\top} x
	+ c_{i} + C>0.
	\end{align*}
	Therefore, the condition set $C_{\delta}$ can be constructed as
	\begin{align*}
	C_{\delta} :=\bigcup_{i = 1}^{N} \{L\in\mathbb{R}^{d_{l}}\mid a_{i}^{\top}L + C >
	\bar{F}^{-1}_{a_{i}^{\top}L}(\delta)\}.
	\end{align*}
	Thus, as the distribution $a_{i}^{\top}L$ is  regularly varying in dimension one for each $i$, we have
	$\limsup_{\delta\rightarrow0}\delta^{-1}\mathrm{P}(L\in C_{\delta})\leq N$, completing the proof.
\end{proofenv}
\end{revisionenv}

\subsection{Proofs for Section \ref{sec-example}}
\begin{revisionenv}
\begin{proofenv}{Proof of Proposition \ref{prop-port-opt}}
Let $\phi(x,l) = \sum_{i=1}^{d} (l_i/x_i) - \eta$ and $\pi(x) = \min_{i=1}^{d} x_i$. The level set is $\Pi = \{x\in\mathbb{R}_{++}^d\mid \min_{i = 1,\ldots, d} x_i = 1\}$. Let $h(\alpha) = \alpha$ and $r(\alpha) = 1$, it follows that $\frac{1}{r(\alpha)}\phi(\alpha\cdot x, h(\alpha)\cdot l) = \phi(x, l)$. In view of the inequalities $
\phi(x,l) \leq \mathbf{1}^{\top} l - \eta$ and $\phi(x,l) \geq \min_{i = 1,\ldots,d} l_i-\eta$ when $x\in\Pi$, we choose the asymptotic uniform bounds as
\begin{align*}
    \Phi_{+}(l) =  \mathbf{1}^{\top} l  - \eta,\qquad\qquad
    \Phi_{-}(l) = \min_{i = 1,\ldots,d} l_i - \eta.
\end{align*}
Furthermore, by definition we construct two approximation sets as
\begin{align*}
    C_{\varepsilon,+} 
    = 
    \left\{l\in\mathbb{R}_{++}^d ~\middle\vert~ \mathbf{1}^{\top} l  \geq \eta -\varepsilon
    \right\},\quad\quad
    C_{\varepsilon,-} 
    = 
    \left\{l\in\mathbb{R}_{++}^d ~\middle\vert~ \min_{i = 1,\ldots, d} l_i \geq \eta +\varepsilon
    \right\}.
\end{align*}

With all the elements that we have already defined, Assumption \ref{assumption-regular-varying} follows directly from the assumption on distribution of $L$. Now we turn to verify Assumption \ref{assumption-level-set}. As $\pi(\alpha\cdot x)=\alpha\cdot\pi(x)$ due to the definition of $\pi(x)$, it suffices to prove that $\lim_{\delta\rightarrow0}\inf_{x\in F_{\delta}}\pi(x)=+\infty$. In view of $\phi(x,L)\leq \mathbf{1}^{\top} L /\pi(x) - \eta$, we have 
\begin{align*}
\label{eq-portfolio-opt-feasible-region}
\begin{array}{ll}
F_\delta = \left\{
x\in\mathbb{R}_{++}^{d}
~
\middle\vert~ \mathrm{P}(\phi(x,L) > 0)\leq \delta
\right\}
&\subseteq
\left\{
x\in\mathbb{R}_{++}^{d}~
\middle\vert~ \mathrm{P}\left(\mathbf{1}^{\top} L> \eta\cdot\pi(x)\right)\leq \delta
\right\}\\
&=\left\{
x\in\mathbb{R}_{++}^{d}~
\middle\vert~ \eta\cdot \pi(x) \geq \bar{F}^{-1}_{\mathbf{1}^{\top}L}(\delta)
\right\}
\end{array}
\end{align*}
Consequently, we have $\inf_{x\in F_{\delta}}\pi(x)\geq \eta^{-1}\bar{F}^{-1}_{\mathbf{1}^{\top}L}(\delta)$. Taking limit for $\delta\rightarrow 0$, we conclude that 
$\lim_{\delta\rightarrow0}\inf_{x\in F_{\delta}}\pi(x)=+\infty$. 

As Assumption \ref{assumption-regular-varying} and \ref{assumption-level-set} are both satisfied, and we also have $\mu(C_{\varepsilon,-}) > 0$, thus Property \ref{condition-outer-approx} is verified due to Lemma \ref{lemma-outer-approx}. In addition, if $\varepsilon\in (0,\eta)$, we have $C_{\varepsilon,+}$ is bounded away from the origin. Thus Assumption \ref{assumption-approx-func} is verified with $S = \mathbb{R}^{d}$.

Finally, we provide closed form expressions for $O_\delta$ and $C_\delta$. Define $\alpha_\delta = \eta^{-1}\cdot \bar{F}^{-1}_{\mathbf{1}^{\top}L}(\delta)$, then it follows that $O_\delta = \bigcup_{\alpha\geq\alpha_{\delta}}\alpha\cdot\Pi = \left\{
x\in\mathbb{R}_{++}^{d}~
\middle\vert~ \eta\cdot \pi(x) \geq \bar{F}^{-1}_{\mathbf{1}^{\top}L}(\delta)
\right\}$, and $C_{\delta} = h(\alpha_{\delta})\cdot\big(C_{\varepsilon, +}\cup K^{c}\cup S^{c}\big) = 
\alpha_{\delta} \cdot C_{\varepsilon, +}
=\left\{l\in\mathbb{R}_{++}^d ~\middle\vert~ \mathbf{1}^{\top} l  \geq (1 -\varepsilon/\eta)\cdot \bar{F}^{-1}_{\mathbf{1}^{\top}L}(\delta)
\right\}.$ By setting $\varepsilon = \eta/2$, we get the expression in the statement of the theorem.
\end{proofenv}

\begin{proofenv}{Proof of Lemma \ref{lemma-minimal-salvage-fund-closed-form}}
We start by showing some properties of $I - Q^{\top}$. Since $Q$ is a non-negative matrix and the row sum is less than 1, it is a sub-stochastic matrix and all of its eigenvalues must be less than 1 in magnitude. This further implies: (1) $I-Q^{\top}$ is invertible; and (2) $(I-Q^{\top})^{-1} = I + \sum_{n = 1}^{\infty} (Q^{\top})^n$ is a non-negative matrix with strictly positive diagonal terms.

Notice that $y = (I-Q^{\top})^{-1}x$ is the unique vector such that $(I-Q^{\top}) y = x$. Let $(y', b')$ be the optimal solution of
\begin{align*}
    \phi\left( x,L\right)  =\min_{y,b}
    \{b \mid(L - y - m) \preceq b\cdot\mathbf{1},\;
    \left(  I-Q^{\top}\right)  y\preceq x, \; y\in\mathbb{R}^{d}_{+}, b\in\mathbb{R}\}
\end{align*}
We have $(I-Q^{\top}) y'\preceq (I-Q^{\top}) y = x$, and we multiply the non-negative matrix $(I-Q^{\top})^{-1}$ on both side, yielding $y'\preceq y$. Obviously, Let $b = \max_{i = 1,\ldots, d} (L_i - y_i)$ such that $(y, b)$ is a feasible solution to above problem. Obviously, it follows from  $y'\preceq y$ that $b' = \max_{i = 1,\ldots, d} (L_i - y'_i - m_i) \geq \max_{i = 1,\ldots, d} (L_i - y_i - m_i) = b$, thus $(y,b)$ is also optimal, which completes the proof. 
\end{proofenv}

\begin{proofenv}{Proof of Proposition \ref{prop-minimal-salvage-fund}}
	Assumption \ref{assumption-regular-varying} follows directly from the
	assumptions of the example. Now we turn to verify Assumption
	\ref{assumption-linear-approx}. Using Lemma \ref{lemma-minimal-salvage-fund-closed-form}, we define $\phi_{-}(x,l) = \phi(x,l) = \max_{i = 1,\ldots, d}\; L_i - \mathbf{e}^{\top}_{i}(I-Q^{\top})^{-1}x - m_i$. Therefore, Assumption \ref{assumption-linear-approx} is satisfied with $N = d$, $a_i = \mathbf{e}_i$, $b_i = -(I-Q)^{-1}\mathbf{e}_i$, $c_i = -m_i$ and $C = 0$. Plugging above values into the expressions of $O_\delta$ and $C_\delta$ given in Theorem \ref{thm-convex}, we get the expressions shown in the statement of the proposition.
\end{proofenv}
\end{revisionenv}

The following lemma is used in the proof of Proposition \ref{prop-quadratic}.
\begin{lemma}
	\label{lemma-geometry}  There exist sets $S_{1},\ldots, S_{2d_{l}}%
	\subseteq\mathbb{R}^{d_{l}}$ with positive Lebesgue measure such that for any
	$z\in\mathbb{R}^{d_{l}}$ with $\|z\|_{2} = 1$, there exist some $S_{i}%
	\subseteq\{l\in\mathbb{R}^{d_{l}}|z^{\top}l> 1\}$.
\end{lemma}

\begin{proofenv}{Proof of Lemma \ref{lemma-geometry}.}
	Let $\mathbf{e}_{i}$ denote the unit vector on the $i$th coordinate in
	$\mathbb{R}^{d_{l}}$ for $i = 1,\ldots, d_{l}$. Fix $z = (z_{1},\ldots,
	z_{d_{l}})\in\mathbb{R}^{d_{l}}$ with $\|z\|_{2} = 1$, define $\theta_{i}$ be
	the angle between $z$ and $\mathbf{e}_{i}$, which satisfies $\cos(\theta_{i})
	= z^{\top}e_{i}$. Since we have $\sum_{i = 1}^{n}\cos(\theta_{i})^{2} = 1$, so
	there exist some $i$ such that $\cos(\theta_{i})^{2} \geq1/n$, thus $z_{i}
	\in[-1,-1/\sqrt{n}] \cup[1/\sqrt{n}, 1]$.  Then, define
	\begin{align*}
	S_{2i - 1}  & = \{l=(l_{1},\ldots,l_{d_{l}})\in\mathbb{R}^{d_{l}}\mid l_{i}>0,
	\;l^{2}_{i} \geq(n-1)\sum_{j\neq i}l_{j}^{2}\},\\
	S_{2i}  & = \{l=(l_{1},\ldots,l_{d_{l}})\in\mathbb{R}^{d_{l}}\mid l_{i}<0,
	\;l^{2}_{i} \geq(n-1)\sum_{j\neq i}l_{j}^{2}\}.
	\end{align*}
	we have either $S_{2i-1} \subset\{l\in\mathbb{R}^{d_{l}}|z^{\top}l> 1\}$ or
	$S_{2i} \subset\{l\in\mathbb{R}^{d_{l}}|z^{\top}l> 1\}$. Thus the proof is complete.
\end{proofenv}
\begin{proofenv}{Proof of Proposition \ref{prop-quadratic}.}
	For the first statement, since $x^{\top}Qx\geq0$ and $A^{\top}x\in\mathbb{R}^{d_{l}%
	}$, and invoking the assumption that $L$ has a positive density,
	\[
	\min_{y\in\mathbb{R}^{d_{l}}\backslash\{\mathbf{0}\}}P(y^{\top} L>0) \geq\min
	_{y:\|y\|_2 = 1}P(y^{\top} L>0)> 0.
	\]
	For the second statement, Assumption
	\ref{assumption-regular-varying} is easy to verify. Notice that $\alpha
	^{-2}\phi(\alpha\cdot x,\alpha\cdot L) = \phi( x, L)$ for all $\alpha>0$, so
	we pick the scaling rate function as $h(\alpha) = \alpha$ and $r(\alpha) =
	\alpha^{2}$. Let $\lambda_{\max}$ denote the maximal eigenvalue of $Q$, and
	$\lambda_{\min}$ denote the minimal eigenvalue of $Q$. The rest of the proof
	will be divided into two cases.
	
	Case 1 ($\lambda_{\max} < 0$): We pick the level set as $\Pi= \{x\in
	\mathbb{R}^{d_{x}}\mid\|x\|_2 = 1 \}$. Since $\lim_{\delta\rightarrow0}
	\inf_{x\in F_{\delta}} \|x\|_2 = \infty$, Assumption \ref{assumption-level-set}
	is verified. Next, we directly show  Property \ref{condition-outer-approx}
	instead of using Lemma \ref{lemma-outer-approx}. For any $x\in\alpha\cdot\Pi$
	we have
	\begin{align*}
	\min_{x\in\alpha\cdot\Pi}\mathrm{P}(x^{\top}Qx +x^{\top}AL > 0)  & \geq\min_{x\in
		\Pi} \mathrm{P}\big(\alpha\lambda_{\min} + x^{\top}AL > 0\big)\\
	&  =\min_{x\in\Pi} \mathrm{P}\bigg(\frac{x^{\top}AL}{\|A^{\top}x\|_{2}} >
	\frac{-\alpha\lambda_{\min}}{\|A^{\top}x\|_{2}}\bigg)\\
	& \geq\min_{z:\|z\| = 1} \mathrm{P}\bigg(z^{\top}L > -\alpha\sigma^{-1}\lambda_{\min
	}\bigg)\\
	(\mbox{Apply Lemma \ref{lemma-geometry}})\quad &  \geq\min_{i = 1,\ldots
		,2d_{l}} \mathrm{P}(L\in-\alpha\sigma^{-1}\lambda_{\min}S_{i} ).
	\end{align*}
	Thus, $\alpha_{\delta}$ can be chosen such that $\alpha_{\delta}= O(\delta)$,
	and $\min_{i = 1,\ldots,2d_{l}} \mathrm{P}(L\in-\alpha_{\delta}\sigma^{-1}\lambda
	_{\min}S_{i} ) > \delta$. As a result, Property \ref{condition-outer-approx}
	is verified. We next turn to derive the asymptotic uniform bound $\Psi_{+}$.
	Observing that
	\begin{align*}
	\sup_{x\in\Pi} \phi(x,L) \leq\lambda_{\max}+\|A\|_{F}\|L\|_2,
	\end{align*}
	we define $\Psi_{+}(L) := \lambda_{\max} +\|A\|_{F}\|L\|_2$. Assumption
	\ref{assumption-approx-func} now follows from the definition of
	$\Psi_{+}$.
	
	Case 2 ($\lambda_{\max} \geq0$): The level set $\Pi$ is chosen as an unbounded
	set $\Pi= \{x\in\mathbb{R}^{d_{x}}\mid x^{\top}Qx = -\|x\|_2\}$ and we have
	$\min_{x\in\Pi} \|x\|_{2} = 1/|\lambda_{\min}|$. For any $x\in\alpha\cdot\Pi$
	we have
	\begin{align*}
	\min_{x\in\alpha\cdot\Pi}\mathrm{P}(x^{\top}Qx +x^{\top}AL > 0)  & \geq\min_{x\in
		\Pi} \mathrm{P}\big(x^{\top}AL > \alpha\big),\\
	&  =\min_{x\in\Pi} \mathrm{P}\bigg(\frac{x^{\top}AL}{\|A^{\top}x\|_{2}} >
	\frac{\alpha}{\|A^{\top}x\|_{2}}\bigg)\\
	&  \geq\min_{z:\|z\| = 1} \mathrm{P}\bigg(z^{\top}L > -\alpha
	\sigma^{-1}\lambda_{\min}\bigg)\\
	(\mbox{Apply Lemma \ref{lemma-geometry}})\quad &  \geq\min_{i = 1,\ldots
		,2d_{l}} \mathrm{P}(L\in-\alpha\sigma^{-1}\lambda_{\min}S_{i} ).
	\end{align*}
	Thus we can pick an $\alpha_{\delta}$ that satisfies Property
	\ref{condition-outer-approx}. Now note that $\sup_{x\in\Pi} \phi(x,L)$ is
	 bounded by
	\begin{align*}
	& \sup_{x\in\Pi} \phi(x,L) \leq\sup_{x\in\Pi} \|x\|_2(\|AL\|_2-1) \leq-\frac{1}%
	{2}|\lambda_{\min}|^{-1} \cdot I(\|AL\|_2 \leq1/2) +\infty\cdot I(\|AL\|_2 > 1)
	\end{align*}
	so we can pick $\Psi_{+}(L) := -\frac{1}{2}|\lambda_{\min}|^{-1} \cdot
	I(\|AL\|_2 \leq1/2) +\infty\cdot I(\|AL\|_2 > 1)$. Consequently Assumption
	\ref{assumption-approx-func} follows immediately.
\end{proofenv}

\clearpage
\section{Additional Numerical Results}
\label{sec-additional-result}
\subsection{Additional Results for Minimal Salvage Fund}
\label{sec-additional-result-minimal-fund}
\begin{figure}[h!]
    \centering
    \begin{subfigure}[b]{0.47\textwidth}
    \centering
    \includegraphics[width=1.1\textwidth]{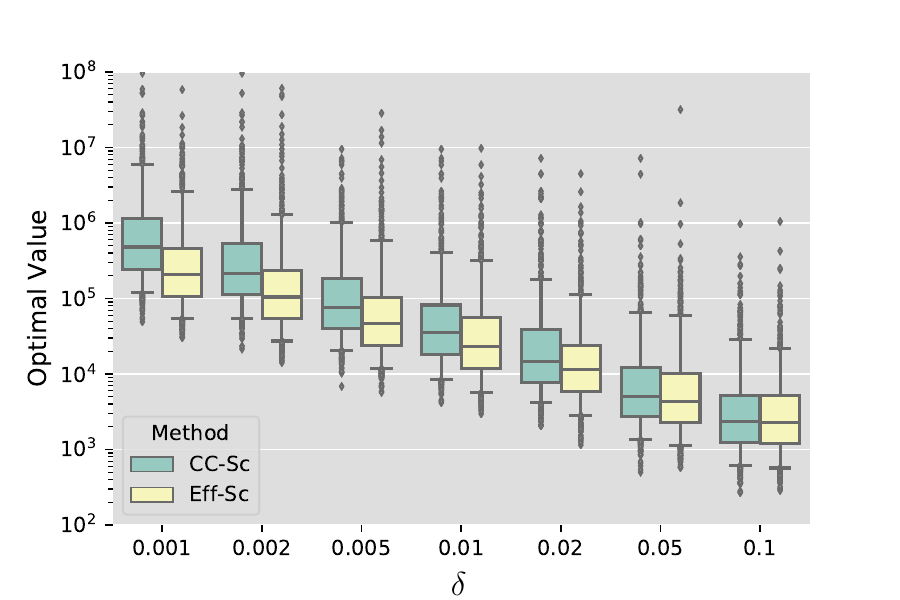}
    \caption{Optimal Value of the Sampled Problem}
    \label{subfig-salvage-fund-obj-val-d5}
    \end{subfigure}
    \hfill
    \begin{subfigure}[b]{0.47\textwidth}
    \centering
    \includegraphics[width=1.1\textwidth]{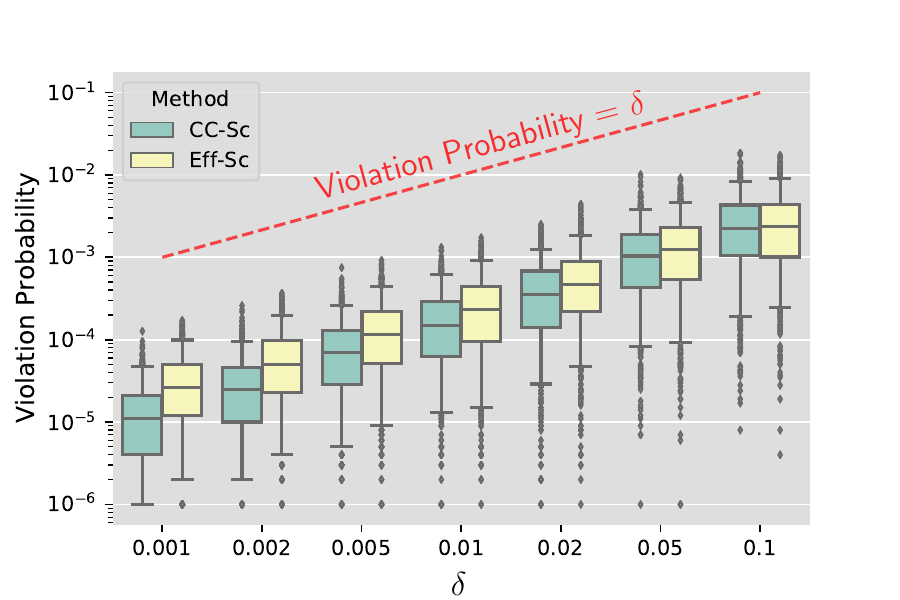}
    \caption{Violation Probability of the Solution}
    \label{subfig-salvage-fund-prob-d5}
    \end{subfigure}
    \caption{Comparison of the quality of optimal solutions given by Eff-Sc and CC-Sc for $d = 5$, in terms of the optimal value shown in Figure \ref{subfig-salvage-fund-obj-val-d5} and the solutions' violation probabilities shown in Figure \ref{subfig-salvage-fund-prob-d5}. 
    }
    \label{fig-solution-quality-salvage-fund-d5}
\end{figure}

\begin{figure}[h!]
    \centering
    \begin{subfigure}[b]{0.47\textwidth}
    \centering
    \includegraphics[width=1.1\textwidth]{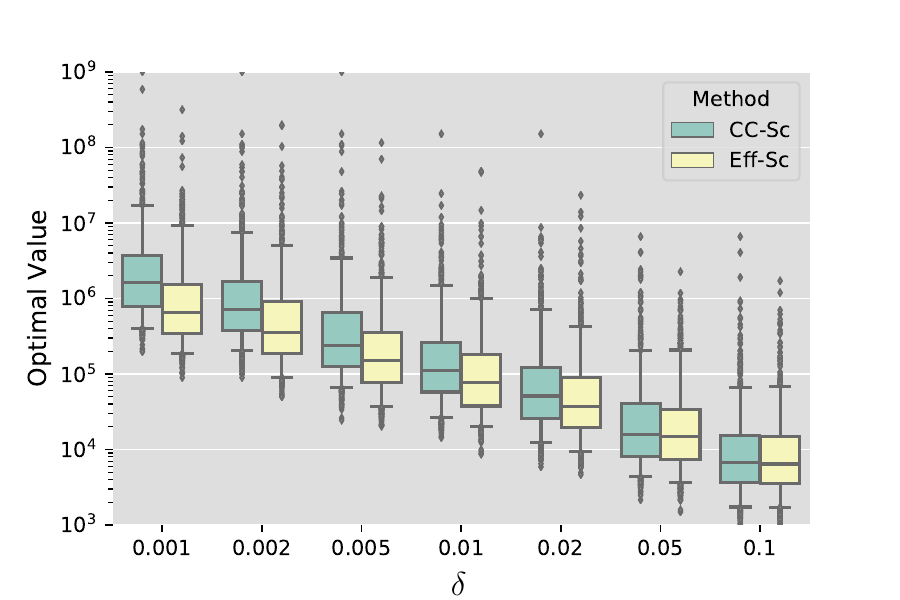}
    \caption{Optimal Value of the Sampled Problem}
    \label{subfig-salvage-fund-obj-val-d10}
    \end{subfigure}
    \hfill
    \begin{subfigure}[b]{0.47\textwidth}
    \centering
    \includegraphics[width=1.1\textwidth]{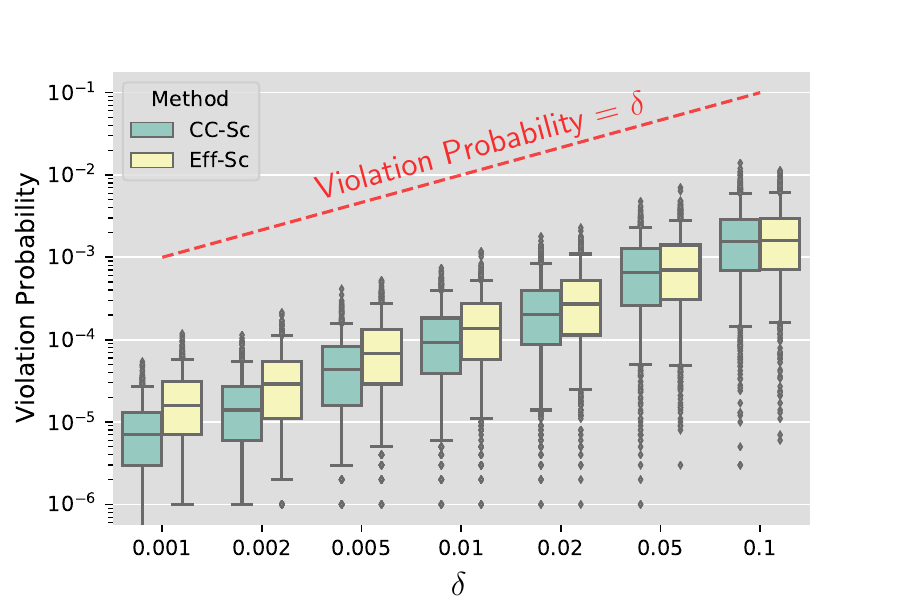}
    \caption{Violation Probability of the Solution}
    \label{subfig-salvage-fund-prob-d10}
    \end{subfigure}
    \caption{Comparison of the quality of optimal solutions given by Eff-Sc and CC-Sc for $d = 10$, in terms of the optimal value shown in Figure \ref{subfig-salvage-fund-obj-val-d10} and the solutions' violation probabilities shown in Figure \ref{subfig-salvage-fund-prob-d10}. 
    }
    \label{fig-solution-quality-salvage-fund-d10}
\end{figure}

In this section we presents additional numerical experiments that demonstrates the quality of the solutions produced by Eff-Sc is better than the solutions produced by CC-Sc. See Figure \ref{fig-solution-quality-salvage-fund-d5} for dimension $d = 5$ and Figure \ref{fig-solution-quality-salvage-fund-d10} for $d = 10$.


\end{document}